\documentclass[11pt,epsf]{article}

\usepackage{comment}
\newcommand {\dfn} {\stackrel{\Delta} {=}}

\newcommand{\eqa}{\stackrel{\mbox{(a)}}{=}}

\newcommand{\geb}{\stackrel{\mbox{(b)}}{\ge}}

\newcommand {\reals} {{\rm I\!R}}

\newcommand {\bq} {\mbox{\boldmath $q$}}
\newcommand {\br} {\mbox{\boldmath $r$}}

\newcommand {\bx} {\mbox{\boldmath $x$}}

\newcommand {\bE} {\mbox{\boldmath $E$}}

\newcommand {\bX} {\mbox{\boldmath $X$}}

\newcommand{\calA}{{\cal A}}

\newcommand{\calC}{{\cal C}}

\newcommand{\calH}{{\cal H}}
\newcommand{\calI}{{\cal I}}

\newcommand{\calS}{{\cal S}}

\newcommand{\calX}{{\cal X}}

\begin{document}
\thispagestyle{empty}
\title{Two New Families of Local Asymptotically Minimax Lower Bounds in
Parameter Estimation}

\author{Neri Merhav
%\thanks{
%Currently on sabbatical leave at HP Laboratories,
%1501 Page Mill Road, MS 3U-4, Palo Alto CA 94304, USA.}
}
\date{}
\maketitle

\begin{center}
The Andrew \& Erna Viterbi Faculty of Electrical Engineering\\
Technion - Israel Institute of Technology \\
Technion City, Haifa 32000, ISRAEL \\
E--mail: {\tt merhav@ee.technion.ac.il}\\
\end{center}
\vspace{1.5\baselineskip}
\setlength{\baselineskip}{1.5\baselineskip}

\abstract{We propose two families of asymptotically local minimax lower
bounds on parameter estimation performance. The first family of bounds applies to any
convex, symmetric loss function that depends solely on the difference between
the estimate and the true underlying parameter value (i.e., the estimation
error), whereas the second is more specifically oriented to the moments of the
estimation error.
The proposed bounds are relatively
easy to calculate numerically (in the sense that their optimization is over relatively few
auxiliary parameters),
yet they turn out to be tighter (sometimes significantly so) than previously reported
bounds that are associated with similar calculation efforts, across a variety
of application examples. In addition to their relative simplicity, they also
have the following advantages: (i) Essentially no regularity conditions are
required regarding the parametric family of distributions; (ii) The bounds are
local (in a sense to be specified); (iii) The bounds provide
the correct order of decay as functions of the number of observations, at
least in all
examples examined; (iv) At least the first family of bounds extends
straightforwardly to vector parameters.}

\section{Introduction}
\label{intro}

The theory of parameter estimation consists of a very large plethora of lower
bounds (as well as upper bounds), that characterize fundamental performance limits of any estimator in a
given parametric model. In this context, it is common to distinguish between
Bayesian bounds (see, e.g., the Bayesian Cram\'er-Rao bound \cite{VanTrees68}, 
the Bayesian Bhattacharyya bound, the Bobrovsky-Zakai bound \cite{BZ76}, the
Belini-Tartara bound \cite{BT74} and the Chazan-Zakai-Ziv bound \cite{CZZ75},
the Weiss-Weinstein bound \cite{Weiss85}, \cite{WW85} and more, see \cite{VTB07}
for a comprehensive overview),
and non-Bayesian bounds,
where in  the former, the parameter to
be estimated is considered a random variable with a given probability
law, as opposed to the latter, where it is assumed an unknown deterministic
constant. The category of non-Bayesian bounds is further subdivided into two
subclasses, one is associated with local bounds that hold for classes of estimators
with certain limitations, such as unbiased estimators (see, e.g., 
the Cram\'er-Rao bound, \cite{Fisher22},
\cite{Dugue37}, \cite{Frechet43}, \cite{Rao45}, \cite{Cramer46}, the
Bhattacharyya bound \cite{Bhattacharyya46}, the Barankin bound
\cite{Barankin49}, the Chapman-Robbins bound
\cite{CR51}, the Fraser-Guttman bound \cite{FG52}, the Keifer bound
\cite{Kiefer52}, and more), and the other is the subclass of minimax bounds
(see, e.g., Ziv and Zakai \cite{ZZ69}, Hajek \cite{Hajek72}, Le Cam
\cite{LeCam73}, Assouad \cite{Assouad83}, Fano \cite{Fano61}, Lehmann \cite[Sections
4.2--4.4]{Lehmann83}, Nazin \cite{Nazin91}, Yang and Barron \cite{YB99},
Guntuboyina \cite{Guntuboyina11a} \cite{Guntuboyina11b},
Kim \cite{Kim20}, and many more).

In this paper, we focus on the minimax approach, and more concretely, on the
local minimax approach. According to the minimax approach, we are
given a parametric family of probability density functions (or probability
mass functions, in the discrete alphabet case),
$\{p(x_1,\ldots,x_n|\theta),~(x_1,\ldots,x_n)\in\reals^n,~\theta\in\Theta\}$,
where $\theta$ is the parameter,
$\Theta$ is the parameter space, 
$n$ is a positive integer designating the number of observations, and we
define a loss
function, $\ell(\theta,\hat{\theta}_n)$, where $\hat{\theta}_n$ is an estimator, which is a
function of the observations $x_1,\ldots,x_n$, only. The minimax performance is defined as
\begin{equation}
R_n(\Theta)\dfn\inf_{\hat{\theta}_n(\cdot)}\sup_{\theta\in\Theta}\bE_\theta\{\ell(\theta,\hat{\theta}_n)\},
\end{equation}
where $\bE_\theta$ denotes expectation w.r.t.\ $p(\cdot|\theta)$. 
As customary, we consider here
loss functions with the property that  $\ell(\theta,\hat{\theta}_n)$ depends on
$\theta$ and $\hat{\theta}_n$ only via their difference, that is,
$\ell(\theta,\hat{\theta}_n)=\rho(\theta-\hat{\theta}_n)$, where in the case of a
scalar parameter, which is the case considered throughout most of this
article, $\rho(\varepsilon)$, $\varepsilon\in\reals$, is
a non-negative function, monotonically non-increasing for $\varepsilon\le 0$, monotonically
non-decreasing for $\varepsilon\ge 0$, and $\rho(0)=0$. 
The local asymptotic minimax performance at the point $\theta\in\Theta$ 
is defined as follows (see also, e.g., \cite{Hajek72}). Let $\{\zeta_n^*,~n\ge 1\}$ be a
positive sequence, tending to infinity, with the property that
\begin{equation}
r(\theta)\dfn\lim_{\delta\downarrow 0}\liminf_{n\to\infty}
\inf_{\hat{\theta}_n(\cdot)}\sup_{\{\theta':~|\theta'-\theta|\le\delta\}}
\zeta_n^*\cdot\bE_{\theta'}\{\rho(\hat{\theta}_n-\theta')\}
\end{equation}
is a strictly positive finite constant. Then, we say that $r(\theta)$ is the
local asymptotic minimax performance with respect to (w.r.t.) $\{\zeta_n^*\}$ at the point
$\theta\in\Theta$. Roughly speaking, the significance is that the performance
of a good estimator, $\hat{\theta}_n$, at $\theta$ is about $R_n(\theta,g_n)\approx
\frac{r(\theta)}{\zeta_n^*}$.
For example, in the mean square error (MSE) case, where
$\rho(\varepsilon)=\varepsilon^2$, and where the observations are Gaussian, i.i.d., with mean $\theta$ and known
variance $\sigma^2$, it is actually shown in \cite[Example 2.4, p.\ 257]{Lehmann83} 
that $r(\theta)=\sigma^2$ w.r.t.\ $\zeta_n^*=n$,
for all $\theta\in\reals$, which is attained by the sample mean estimator,
$\hat{\theta}_n=\frac{1}{n}\sum_{i=1}^nx_i$.

Our focus on this work is on the derivation of some new lower bounds that are:
(i) essentially free of regularity conditions on the smoothness of
the parametric family $\{p(\cdot|\theta),~\theta\in\Theta\}$, (ii)
relatively simple and easy to calculate, at least numerically, which amounts
to the property that the bound contains only a small number of auxiliary parameters 
to be numerically optimized (typically, no more than two or three
parameters), (iii) tighter than earlier reported bounds that are
associated with similar calculation efforts as described in (ii), and (iv) lend themselves to
extensions that yield even stronger bounds (albeit with more auxiliary
parameters to be optimized), as well as extensions to
vector parameters. We propose two families of lower
bounds on $R_n(\Theta)$, along with their local versions, of bounding $r(\theta)$, 
with the four above described properties. The first applies to any
convex, symmetric loss function $\rho(\cdot)$,
whereas the second is more specifically oriented to the moments of the
estimation error, $\rho(\varepsilon)=|\varepsilon|^t$, where $t$ is a positive real, not
necessarily an integer, with special attention devoted to the MSE case, $t=2$.

To put this work in the perspective of earlier work on minimax estimation, we
next briefly review some of the basic approaches in this problem area. Admittedly,
considering the vast amount of literature on the subject,
our review below is by no means exhaustive. For a more
comprehensive review, the reader is referred to Kim \cite{Kim20}. 

First, observe the simple fact that
the minimax performance is lower bounded by the Bayesian
performance of the same loss function (see, e.g., 
\cite{VanTrees68}, \cite{BZ76},
\cite{BT74}, \cite{CZZ75},
\cite{Weiss85}, \cite{WW85}, \cite{VTB07})
for any prior on the parameter,
$\theta$, and so, every lower bound on the Bayesian performance is
automatically a valid lower bound also on the minimax performance
\cite[Section 4.2]{Lehmann83}. Indeed, in \cite[Section 2.3]{Guntuboyina11a}
it is argued that the vast majority of existing minimax lower bounding techniques are
based upon bounding the Bayes risk from below w.r.t.\ some prior. 
Many of these Bayesian bounds, however, are subjected to certain restrictions
and regularity conditions concerning the smoothness of the prior and the family
of densities, $\{p(\cdot|\theta),~\theta\in\Theta\}$.

Dating back
to Ziv and Zakai's 1969 article \cite{ZZ69} on parameter estimation, applied
mostly in the context of time-delay estimation, this prior puts all its mass
equally on two values, $\theta_0$ and $\theta_1$, of the parameter $\theta$,
and considering an hypothesis testing problem of distinguishing between the two
hypotheses, $\calH_0:~\theta=\theta_0$ and
$\calH_1:~\theta=\theta_1$ with equal priors. A simple argument regarding the
sub-optimality of a decision rule that is based on estimating $\theta$ and
deciding on the hypothesis with the closer value of $\theta$, combined with
Chebychev's inequality, yields a simple lower bound on the corresponding Bayes risk,
and hence also the minimax risk, in terms of the probability of error of the
optimal decision rule. Five years later, Bellini and Tartara \cite{BT74}, and
then independently, Chazan, Zakai and Ziv \cite{CZZ75}, improved the bound of
\cite{ZZ69} using somewhat different arguments, and obtained Bayesian bounds
that apply to the uniform prior. These bounds are also given in terms of the
error probability pertaining to the optimal MAP decision rule of binary hypothesis
testing with equal priors, but this
time, it had an integral form. These bounds were demonstrated to be fairly
tight in several application examples, but they suffer from two main drawbacks:
(i) they are difficult to calculate in most cases, (ii) it is not apparent how
to extend these bounds to apply to general, non-uniform priors.
Shortly before the Bellini-Tartara and
the Chazan-Zakai-Ziv articles were published, Le Cam \cite{LeCam73} proposed a minimax
lower bound, which is also given in terms of the error probability associated
with binary hypothesis testing, or equivalently, the total variation between
$P(\cdot|\theta_0)$ and $P(\cdot|\theta_1)$, under the postulate that the loss
function $\ell(\cdot,\cdot)$ is a metric. We will refer to Le Cam's bound
in a more detailed manner later, in the context of our first proposed bound.
A decade later, Assouad
\cite{Assouad83} extended Le Cam's two-points testing bound to multiple points
and devised the so called hypercube method. Another, related bounding
technique, that is based on multiple
test points, and referred to as Fano's method, amounts to further bounding
from below the
error probability of multiple hypotheses using Fano's inequality
\cite[Sect.\ 2.10]{CT06}. Considering the large number of auxiliary parameters
to be optimized when multiple hypotheses are present, these bounds demand
heavy computational efforts. Also, Fano's inequality is often loose, even
though it is adequate enough for the purpose of proving converse-to-coding theorems
in information theory \cite{CT06}. In later years, Le Cam \cite{LeCam86} and Yu
\cite{Yu97} extended Le Cam's original approach to apply to testing mixtures
of densities. More recently, Yang and Barron \cite{YB99}
related the minimax problem to the metric entropy of the parametric family,
$\{p(\cdot|\theta_,~\theta\in\Theta\}$, and
Cai and Zhou \cite{CZ12} combined Le Cam's
and Assouad's methods by considering a larger number of dimensions. 
Guntuboyina \cite{Guntuboyina11a}, \cite{Guntuboyina11b} pursued a different
direction by deriving minimax lower bounds using
$f$-divergences.

The outline of this article is as follows.
In Section 2, we define the problem setting, provide a few formal definitions
and establish the notation. In Section 3, we develop the first family of bounds
and finally, in Section 4, we present the second family.

\section{Problem Setting, Definitions and Notation}
\label{psdn}

Consider a family of probability density functions (pdf's),
$\{p(\cdot|\theta),~\theta\in\Theta\}$,
where $\theta$ is a scalar parameter to be
estimated and $\Theta\subseteq\reals$ is the parameter space.
We denote by $\bE_\theta\{\cdot\}$ the expectation operator
w.r.t.\ $p(\cdot|\theta)$.
Let $\bX=(X_1,\ldots,X_n)$ be a random vector
of observations,
governed by $p(\cdot|\theta)$ for some $\theta\in\Theta$. The support of
$p(\cdot|\theta)$ is assumed
$\calX^n$, the $n$th Cartesian power of the alphabet, $\calX$ of each
component, $X_i$, $i=1,\ldots,n$. The alphabet $\calX$ may be a finite set, a countable set, a
finite interval, an infinite interval, or the entire real line. In the first
two cases, the pdfs should be understood to replaced by probability
mass functions and integrations over the observation space should be replaced
by summations.
A realization of $\bX$ will be denoted by $\bx=(x_1,\ldots,x_n)\in\calX^n$.

An estimator, $\hat{\theta}_n$, is given by any function of the observation vector,
$g_n:\calX^n\to\Theta$, that is, $\hat{\theta}_n=g_n[\bx]$. Since $\bX$ is
random, then so is the estimate, $g_n[\bX]$, as well as the estimation error,
$\varepsilon_n[\bX]=\hat{\theta}_n-\theta=g_n[\bX]-\theta$. We associate with
every value possible value, $\epsilon$, of $\varepsilon_n[\bX]$ a certain loss (or `cost', or `price'),
$\rho(\epsilon)$, where $\rho(\cdot)$ is a non-negative function with the
following properties: (i) monotonically non-increasing for $\epsilon\le 0$,
(ii) monotonically non-decreasing for $\epsilon\ge 0$, and (iii) $\rho(0)=0$.

In Section \ref{convexsymmetric}, we assume, in addition, that $\rho(\cdot)$ is: (iv) convex, and
(v) symmetric, i.e., $\rho(-\epsilon)=\rho(\epsilon)$ for every $\epsilon$.
In Section \ref{moments}, we assume, more specifically, that
$\rho(\varepsilon)=|\varepsilon|^t$,
where $t$ is a positive constant, not necessarily an integer. This is a special
case of the class of loss functions considered in Section \ref{convexsymmetric}, except when
$t\in(0,1)$, in which case, $|\varepsilon|^t$ is a concave (rather than a convex) function
of $\varepsilon$. 

The expected cost of an estimator $g_n$ at a point $\theta\in\Theta$, is defined as
\begin{equation}
R_n(\theta,g_n)\dfn\bE_\theta\{\rho(g_n[\bX]-\theta)\}.
\end{equation}
The global minimax performance is defined as
\begin{equation}
R_n(\Theta)\dfn\inf_{g_n}\sup_{\theta\in\Theta}R_n(\theta,g_n).
\end{equation}
Another, related notion is that of local asymptotic minimax performance,
defined in Section \ref{intro}, and repeated here for the sake of
completeness.
Let $\{\zeta_n^*,~n\ge 1\}$ be a
positive sequence, tending to infinity, with the property that
\begin{equation}
r(\theta)\dfn\lim_{\delta\downarrow 0}\liminf_{n\to\infty}
\inf_{g_n}\sup_{\{\theta':~|\theta'-\theta|\le\delta\}}\zeta_n^*\cdot\bE_{\theta'}\{\rho(g_n[\bX]-\theta')\}
\end{equation}
is a strictly positive finite constant. Then, we say that $r(\theta)$ is the
local asymptotic minimax performance w.r.t.\ $\{\zeta_n^*\}$ at
$\theta\in\Theta$. The sequence $\{1/\zeta_n^*\}$ is referred to as the
convergence rate of the minimax estimator.

As described in Section \ref{intro}, our objective in this work is to derive relatively
simple and easily computable lower bounds to $r(\theta)$, which are as tight
as possible. While many existing lower bounds in the literature are satisfactory
in terms of yielding the correct rate of convergence, $1/\zeta_n^*$, here we wish
to improve the bound on the constant factor, $r(\theta)$. Many of our examples involve
numerical calculations which include optimization over auxiliary parameters
and occasionally also
numerical integrations. All these calculations were carried out using MATLAB.

\section{Lower Bounds for Convex Symmetric Loss Functions}
\label{convexsymmetric}

We begin from the following simple bound.\\

\noindent
{\bf Theorem 1.} 
Let the assumptions of Section \ref{psdn} be satisfied and
let $\rho(\cdot)$ be a symmetric convex loss function.
Then,
\begin{equation}
\label{1stlowerbound}
R_n(\Theta)\ge
\sup_{\theta_0,~\theta_1\in\Theta}\left\{2\cdot\rho\left(\frac{\theta_1-\theta_0}{2}\right)\cdot
\sup_{0\le q\le 1}P_{\mbox{\tiny
e}}(q,\theta_0,\theta_1)\right\},
\end{equation}
where $P_{\mbox{\tiny e}}(q,\theta_0,\theta_1)$ is defined as
\begin{equation}
P_{\mbox{\tiny e}}(q,\theta_0,\theta_1)\dfn\int_{\reals^n}\min\{q\cdot
p(\bx|\theta_0),(1-q)\cdot p(\bx|\theta_1)\}\mbox{d}\bx
\end{equation}
which is identified as the error probability
associated with the optimal maximum a posteriori (MAP) decision rule for the binary hypothesis
testing problem problem, $\calH_0:~\bX\sim p(\cdot|\theta_0)$ and 
$\calH_1:~\bX\sim p(\cdot|\theta_1)$ with priors $q$ and $1-q$,
respectively.\\

\noindent
{\em Proof of Theorem 1.}
For every $\theta_0,\theta_1\in\Theta$ and $q\in[0,1]$,
\begin{eqnarray}
R_n(\Theta)&\ge&
q\bE_{\theta_0}\{\rho(g_n[\bX]-\theta_0)\}+
(1-q)\bE_{\theta_1}\{\rho(g_n[\bX]-\theta_1)\}\nonumber\\
&\eqa&
\int_{\reals^n}
\bigg[q\cdot p(\bx|\theta_0)\rho(g_n[\bx]-\theta_0)+\nonumber\\
& &(1-q)\cdot p(\bx|\theta_1)\rho(\theta_1-g_n[\bx])\bigg]
\mbox{d}\bx\nonumber\\
&\ge&
2\cdot\int_{\reals^n}
\min\{q\cdot p(\bx|\theta_0),(1-q)\cdot p(\bx|\theta_1)\}\times\nonumber\\
& &\left[\frac{1}{2}\rho(g_n[\bx]-\theta_0)+
\frac{1}{2}\rho(\theta_1-g_n[\bx])\right]
\mbox{d}\bx\nonumber\\
&\geb& 2\cdot
\int_{\reals^n}
\min\{q\cdot p(\bx|\theta_0),(1-q)\cdot p(\bx|\theta_1)\}\times\nonumber\\
& &\rho\left(\frac{g_n[\bx]-\theta_0}{2}+
\frac{\theta_1-g_n[\bx]}{2}\right)
\mbox{d}\bx\nonumber\\
&=&2\cdot\rho\left(\frac{\theta_1-\theta_0}{2}\right)\cdot
\int_{\reals^n}
\min\{q\cdot p(\bx|\theta_0),(1-q)\cdot p(\bx|\theta_1)\}\mbox{d}\bx\nonumber\\
&=&2\cdot\rho\left(\frac{\theta_1-\theta_0}{2}\right)\cdot P_{\mbox{\tiny
e}}(q,\theta_0,\theta_1),
\end{eqnarray}
where (a) is due to the assumed symmetry of $\rho(\cdot)$ and
(b) is by its assumed convexity. Since the inequality,
\begin{equation}
R_n(\Theta)\ge 
2\cdot\rho\left(\frac{\theta_1-\theta_0}{2}\right)\cdot P_{\mbox{\tiny
e}}(q,\theta_0,\theta_1)
\end{equation}
applies to every $\theta_0,\theta_1\in\Theta$ and $q\in[0,1]$, it applies, in
particular, also to the supremum over these auxiliary parameters.
This completes the proof of Theorem 1.\\

Before we proceed, a few comments are in order.\\

\noindent
1. Note that $P_{\mbox{\tiny e}}(q,\theta_0,\theta_1)$ is a concave
function of $q$ for fixed $(\theta_0,\theta_1)$, as it
can be presented as the
minimum among a family of affine functions of $q$,
given by $\min_{\Omega}[q\cdot P(\Omega|\theta_1)+(1-q)\cdot
P(\Omega^c|\theta_2)]$, where $\Omega$ runs over all possible subsets of the
observation space, $\calX^n$. Another way to see why this is true is by
observing that it is defined by an
integral, whose integrand, $\min\{q\cdot p(\bx|\theta_0),(1-q)\cdot
p(\bx|\theta_1)\}$, is concave in $q$. Clearly, $P_{\mbox{\tiny
e}}(q,\theta_0,\theta_1)=0$ whenever $q=0$ or $q=1$. Thus, $P_{\mbox{\tiny e}}(q,\theta_0,\theta_1)$
is maximized by some $q$ between $0$ and $1$. If $P_{\mbox{\tiny e}}(q,\theta_0,\theta_1)$
is strictly concave in $q$, then the maximizing $q$ is unique.\\

\noindent 
2. Although we have scalar parameters in mind (throughout 
most of this work), the above proof continues to hold as is also
when $\theta$ is a vector and $\rho$ is symmetric and jointly convex in all components of
the estimation error.\\

\noindent
3. Note that the lower bound (\ref{1stlowerbound}) is tighter than the lower bound of
$\rho(\frac{\theta_1-\theta_0}{2})P_{\mbox{\tiny
e}}(\frac{1}{2},\theta_0,\theta_1)$, that was obtained 
in \cite[eqs.\ (6)-(9a)]{ZZ69}, both because of the factor of 2 and
because of the freedom to optimize $q$ rather than setting $q=1/2$.
In a further development of \cite{ZZ69} the factor of 2 was accomplished too,
but at the price of assuming that
the density of the estimation error is symmetric about the origin (see
discussion after (10) therein), which limits the class of estimators to which
the bound applies. The factor of 2 and the degree of freedom $q$ 
are also the two ingredients that make the difference between
(\ref{1stlowerbound}) and the lower bound
due to Le Cam \cite{LeCam73} (see also \cite{Kim20} and
\cite{Guntuboyina11a}). In \cite[Chapter 2]{Guntuboyina11a} Guntuboyina reviews standard
bounding techniques, including those of Le Cam, Assouad and Fano. In
particular, in Example 2.3.2 therein, Guntiboyina
presents a lower bound in terms of the error probability associated 
with general priors. However, the coefficient
in front of the error probability factor there is given by $\frac{\eta}{2}$, where
in the case of two hypotheses,
$\eta=\min_\vartheta\{\rho(\theta_0-\vartheta)+\rho(\theta_1-\vartheta)\}$,
in our notation.
Now, if $\rho$ is symmetric and monotonically non-decreasing in the absolute
error, then the minimizing $\vartheta$ is given by
$\frac{\theta_0+\theta_1}{2}$, which yields
$\frac{\eta}{2}=\rho\left(\frac{\theta_1-\theta_0}{2}\right)$ and so, again, the
resulting bound is of the same form as (\ref{1stlowerbound}) except 
that it lacks the prefactor of 2.

Our first example demonstrates Theorem 1 on a somewhat technical, but simple
model, with an emphasis on the point that the optimal $q$ may differ from
$1/2$ and that it is therefore useful to maximize
w.r.t. $q$ in order to improve the bound relative to the choice $q=1/2$.\\

\noindent
{\bf Example 1.} 
Let $X$ a random variable distributed exponentially according to
\begin{equation}
p(x|\theta)=\theta e^{-\theta x},~~~x\ge 0,
\end{equation}
and $\Theta=\{1,2\}$, so that the only possibility to select two different
values of $\theta$ in the lower bound are $\theta_0=1$ and $\theta_1=2$.
In terms of the hypothesis testing problem pertaining to the lower bound,
the likelihood ratio test (LRT) is by comparison of $qe^{-x}$ to
$(1-q)\cdot 2 e^{-2x}$.
Now, if $2(1-q) \le q$, or equivalently,
$q\ge\frac{2}{3}$, the decision is always in favor of
$\calH_0$, and then $P_{\mbox{\tiny e}}(q,1,2)=1-q$.
For $q < \frac{2}{3}$, the optimal LRT compares $X$ to $x_0(q)=
\ln\frac{2(1-q)}{q}$. If $X>x_0(q)$, one decides in favor of $\calH_0$,
otherwise -- in favor of $\calH_1$. Thus,
\begin{eqnarray}
P_{\mbox{\tiny
e}}(q,1,2)&=&q\int_0^{x_0(q)}e^{-x}\mbox{d}x+(1-q)\int_{x_0(q)}^\infty
2e^{-2x}\mbox{d}x\nonumber\\
&=&q[1-e^{-x_0(q)}]+(1-q)e^{-2x_0(q)}\nonumber\\
&=&q\left[1-\frac{q}{2(1-q)}\right]+(1-q)\left[\frac{q}{2(1-q)}\right]^2\nonumber\\
&=&\frac{q(4-5q)}{4(1-q)}.
\end{eqnarray}
In summary,
\begin{equation}
P_{\mbox{\tiny
e}}\left(q,1,2\right)=\left\{\begin{array}{ll}
\frac{q(4-5q)}{4(1-q)} & 0\le q < \frac{2}{3}\\
1-q & \frac{2}{3}\le q \le 1\end{array}\right.
\end{equation}
It turns out that for
$q=\frac{1}{2}$,
$P_{\mbox{\tiny e}}(\frac{1}{2},1,2)=\frac{3}{8}=0.375$, whereas the
maximum is $0.382$, attained at $q=1-\frac{\sqrt{5}}{5}=0.5528$.
Thus,
\begin{equation}
R_n(\Theta)\ge 2\cdot 0.382
\cdot\rho\left(\frac{2-1}{2}\right)=0.764\cdot\rho\left(\frac{1}{2}\right).
\end{equation}
This concludes Example 1.\\

In the above example, we considered just one observation, $n=1$. From now on,
we will refer to the case where the $n\gg 1$. In particular,
the following simple corollary to Theorem 1 yields a local asymptotic minimax lower
bound.\\

\noindent
{\bf Corollary 1.}
For a given $\theta\in\Theta$ and a constant $s$, let $\{\xi_n\}_{n\ge 1}$
denote a sequence tending to zero with the
property that 
\begin{equation}
\lim_{n\to\infty}\max_qP_{\mbox{\tiny e}}(q,\theta,\theta+2s\xi_n) 
\end{equation}
exists and is given
by a strictly positive constant, which will be denoted by $P_{\mbox{\tiny
e}}^\infty(\theta,s)$. Also, let 
\begin{equation}
\omega(s)\dfn\lim_{u\to 0}\frac{\rho(s\cdot
u)}{\rho(u)}.
\end{equation}
Then, the local asymptotic minimax performance w.r.t.\ $\zeta_n=1/\rho(\xi_n)$
is lower bounded by
\begin{equation}
r(\theta)\ge \sup_{s\in\reals}\left\{2\omega(s)\cdot P_{\mbox{\tiny
e}}^\infty(\theta,s)\right\}.
\end{equation}
\\

Corollary 1 is readily obtained from Theorem 1 by substituting $\theta_0=\theta$
and $\theta_1=\theta+2s\xi_n$ in eq.\ (\ref{1stlowerbound}), 
then multiplying both sides of the inequality by $\zeta_n=1/\rho(\xi_n)$ 
and finally, taking the limit inferior of both sides.

\subsection{Examples for Corollary 1}

We next study a few examples for the use of Corollary 1. Similarly as in Example 1,
we emphasize again in Example 2 below the
importance of having the degree of freedom to maximize over the prior $q$
rather than to fix $q=\frac{1}{2}$. Also, in all the examples that were examined,
the rate of convergence, $1/\zeta_n=\rho(\xi_n)$ is the same as the optimal
rate of convergence, $1/\zeta_n^*$. 
In other words, it is tight in the sense that
there exists an estimator (for example, the maximum likelihood estimator) 
for which $R_n(\theta,g_n)$ tends to zero at the same rate. In some of these
examples, we compare our lower bound to $r(\theta)$ to those of earlier
reported results on the same models.\\

\noindent
{\bf Example 2.}
Let $X_1,\ldots,X_n$ be independently, identically distributed (i.i.d.) random
variables, uniformly distributed in the range $[0,\theta]$.
In the corresponding hypothesis testing problem of Theorem 1, the hypotheses
are $\theta=\theta_0$ and $\theta=\theta_1 > \theta_0$ with priors, $q$ and $1-q$.
There are
two cases: If $q/\theta_0^n < (1-q)/\theta_1^n$, or equivalently, $q
<\frac{\theta_0^n}{\theta_0^n+\theta_1^n}$,
one decides always in favor of
$\calH_1$ and so, the probability of error is $q$. If, on the other hand,
$q/\theta_0^n >
(1-q)/\theta_1^n$, namely, $q> \frac{\theta_0^n}{\theta_0^n+\theta_1^n}$, we decide
in favor of $\calH_1$ whenever $\max_iX_i>\theta_0$ and
then an error occurs only if $\calH_1$ is true, yet $\max_iX_i< \theta_0$, which
happens with probability $(1-q)\left(\frac{\theta_0}{\theta_1}\right)^n$. Thus,
\begin{equation}
P_{\mbox{\tiny e}}(q,\theta_0,\theta_1)=\left\{\begin{array}{ll}
q & q<\frac{\theta_0^n}{\theta_0^n+\theta_1^n}\\
(1-q)\left(\frac{\theta_0}{\theta_1}\right)^n & q \ge
\frac{\theta_0^n}{\theta_0^n+\theta_1^n}\end{array}\right.=
\min\left\{q,(1-q)\left(\frac{\theta_0}{\theta_1}\right)^n\right\},
\end{equation}
which is readily seen to be maximized by
$q=\frac{\theta_0^n}{\theta_0^n+\theta_1^n}$ and then
\begin{equation}
\max_qP_{\mbox{\tiny
e}}(q,\theta_0,\theta_1)=\frac{\theta_0^n}{\theta_0^n+\theta_1^n}.
\end{equation}
Now, to apply Corollary 1, let $\theta_0=\theta$ and
$\theta_1=\theta_0(1+2\sigma/n)$, which amounts to
$s=\theta_0\sigma=\theta\sigma$ and $\xi_n=1/n$.
Then,
\begin{equation}
P_{\mbox{\tiny
e}}^\infty(\theta,s)=\lim_{n\to\infty}\frac{1}{1+[1+2s/(\theta
n)]^n}=\frac{1}{1+e^{2s/\theta}}.
\end{equation}
In case of the MSE criterion, $\rho(\varepsilon)=\varepsilon^2$, we have $\omega(s)=s^2$, and so,
\begin{equation}
r(\theta)\ge\sup_{s\ge 0}\frac{2s^2}{1+e^{2s/\theta}}=\theta^2\cdot\sup_{u\ge
0}\frac{u^2}{2(1+e^u)}=0.2414\theta^2
\end{equation}
w.r.t.\ $\zeta_n=1/\rho(\xi_n)=1/(1/n)^2=n^2$.
This bound will be further improved 
upon in Section 4.

If instead of maximizing w.r.t.\ $q$, we select $q=1/2$, then
\begin{equation}
P_{\mbox{\tiny
e}}\left(\frac{1}{2},\theta,\theta\left(1+\frac{2\sigma}{n}\right)\right)=
\frac{1}{2}\cdot\left[\frac{\theta}{\theta(1+2\sigma/n)}\right]^n\to
\frac{1}{2}\cdot e^{-2\sigma}=\frac{1}{2}\cdot e^{-2s/\theta},
\end{equation}
and then the resulting bound would become
\begin{equation}
r(\theta)\ge \sup_{s\ge 0}s^2e^{-2s/\theta}=\theta^2\cdot \sup_{u\ge
0}\frac{u^2e^{-u}}{4}=0.1353\theta^2
\end{equation}
w.r.t.\ $\zeta_n=n^2$. Therefore, the maximization over $q$ plays an important
role here in terms of tightening the lower bound to $r(\theta)$.

More generally,
for $\rho(\varepsilon)=|\varepsilon|^t$ ($t\ge 1$), $\omega(s)=|s|^t$
and we obtain
\begin{equation}
r(\theta)\ge \theta^t\cdot\sup_{u\ge 0}\frac{2u^t}{1+e^{2u}},
\end{equation}
w.r.t.\ $\zeta_n=n^t$,
where the supremum, which is in fact a maximum, can always be calculated
numerically for every given $t$.
For large $t$, the maximizing $u$ is approximately $t/2$, which yields 
\begin{equation}
r(\theta)\ge \frac{(t\theta)^t}{2^{t-1}(1+e^t)}.
\end{equation}
On the other hand, for $q=1/2$, we end up with
\begin{equation}
r(\theta)\ge \sup_{s\ge 0}s^te^{-2s/\theta}=\left(\frac{t\theta}{2e}\right)^t.
\end{equation}
For large $t$, the bound of $q=1/2$ is inferior to the bound with the optimal
$q$, by a factor of about $1/2$.\\

\noindent
{\bf Example 3.}
Let $X_1,X_2,\ldots,X_n$ be i.i.d.\ random variables, uniformly distributed in the interval
$[\theta,\theta+1]$. For the hypothesis testing problem, let $\theta_1$ be chosen between $\theta_0$ and
$\theta_0+1$. Clearly, if $\min_iX_i <\theta_1$, the underlying hypothesis is
certainly $\calH_1$. Likewise, if $\max_iX_i > \theta_0+1$, the decision is in
favor of $\calH_0$ with certainty. Thus, an error can occur only if all
$\{X_i\}$ fall in the interval $[\theta_1,\theta_0+1]$, an event that occurs
with probability $(\theta_0+1-\theta_1)^n$. In this event, the best
to be done is to select the hypothesis with the larger prior with a probability
of error given by $\min\{q,1-q\}$. 
Thus,
\begin{equation}
P_{\mbox{\tiny
e}}\left(q,\theta_0,\theta_1\right)=(\theta_0+1-\theta_1)^n\cdot\min\{q,1-q\},
\end{equation}
and so,
\begin{equation}
\max_qP_{\mbox{\tiny
e}}\left(q,\theta_0,\theta_1\right)=\frac{1}{2}[1-(\theta_1-\theta_0)]^n,
\end{equation}
achieved by $q=1/2$. 
Now, let us select $\xi_n=1/n$, which yields
\begin{equation}
\lim_{n\to\infty}\max_qP_{\mbox{\tiny
e}}\left(q,\theta_0,\theta_1\right)=\lim_{n\to\infty}\frac{1}{2}\cdot\left(1-\frac{2s}{n}\right)^n=\frac{1}{2}\cdot
e^{-2s}.
\end{equation}
For $\rho(\varepsilon)=|\epsilon|^t$, ($t\ge 1$), we have
\begin{eqnarray}
r(\theta)\ge \sup_{s\ge 0}2s^t\cdot\frac{1}{2}e^{-2s}=
\sup_{s\ge 0}s^te^{-2s}=\left(\frac{t}{2e}\right)^t
\end{eqnarray}
w.r.t.\ $\zeta_n=1/|1/n|^t=n^t$.
For the case of MSE, $t=2$, $r(\theta)\ge e^{-2}=0.1353$.
The constant $0.1353$ should be compared with 
$\frac{1-1/\sqrt{2}}{128}= 0.0023$ \cite[Example 4.9]{Parag22}, which is two orders of
magnitude smaller. This concludes Example 3.\\

\noindent
{\bf Example 4.}
Let $X_i=\theta+Z_i$, where $\{Z_i\}$ are i.i.d., Gaussian random variables with zero mean
and variance $\sigma^2$. Here, for the corresponding binary hypothesis testing
problem, the optimal value of $q$ is always $q^*=\frac{1}{2}$. This can be
readily seen from the concavity of $P_{\mbox{\tiny e}}(q,\theta_0,\theta_1)$
in $q$ and its symmetry around $q=1/2$, as $P_{\mbox{\tiny
e}}(q,\theta_0,\theta_1)=P_{\mbox{\tiny e}}(1-q,\theta_0,\theta_1)$.
Since
\begin{eqnarray}
P_{\mbox{\tiny
e}}\left(\frac{1}{2},\theta_0,\theta_1\right)&=&\mbox{Pr}\left\{\sum_{i=1}^nZ_i\ge
\frac{n(\theta_1-\theta_0)}{2}\right\}\nonumber\\
&=&\mbox{Pr}\left\{\frac{1}{\sigma\sqrt{n}}\sum_{i=1}^nZ_i\ge\frac{\sqrt{n}(\theta_1-\theta_0)}{2\sigma
}\right\}\nonumber\\
&=&Q\left(\frac{\sqrt{n}(\theta_1-\theta_0)}{2\sigma}\right),
\end{eqnarray}
where
\begin{equation}
Q(t)\dfn \int_t^\infty\frac{e^{-u^2/2}\mbox{d}u}{\sqrt{2\pi}},
\end{equation}
we select $\xi_n=\frac{1}{\sqrt{n}}$, which yields
$\theta_1-\theta_0=\frac{2s}{\sqrt{n}}$ 
\begin{equation}
P_{\mbox{\tiny e}}^\infty(\theta,s)=Q\left(\frac{s}{\sigma}\right),
\end{equation}
and then for the MSE case, $\omega(s)=s^2$, 
\begin{equation}
r(\theta)\ge \sup_{s\ge 0} \left\{2s^2Q\left(\frac{s}{\sigma}\right)\right\}=
\sigma^2\cdot\sup_{u\ge 0}\{2u^2Q(u)\}
=0.3314\sigma^2
\end{equation}
w.r.t.\ $\zeta_n=1/(1/\sqrt{n})^2=n$, and so, the asymptotic lower bound to
$R_n(\theta,g_n)$ is $0.3314\sigma^2/n$.

We now compare this bound (which will be further improved in Section 4) with a few earlier reported results.
In one of the versions of Le Cam's bound \cite[Example
4.7]{Parag22} for the same model, the lower bound to $r(\theta)$ turns
out to be $\frac{\sigma^2}{24}=0.0417\sigma^2$, namely, an order of magnitude smaller. Also,
in \cite[Example 3.1]{Kim20},
another version of Le Cam's method yields
$r(\theta)\ge(1-\sqrt{1/2})\sigma^2/8=0.0366\sigma^2$.
According to \cite[Corollary 4.3]{Rogers18},
$r(\theta)\ge\frac{\sigma^2}{8e}=0.046\sigma^2$.
Yet another comparison is with \cite[Theorem 5.9]{Rigollet15}, 
where we find an inequality, which in our notation reads as follows:
\begin{equation}
\sup_{\theta\in\Theta}P_\theta\left\{(g_n[\bX]-\theta)^2\ge
\frac{2\alpha\sigma^2}{n}\right\}\ge
\frac{1}{2}-\alpha~~~~\alpha\in\left(0,\frac{1}{2}\right).
\end{equation}
Combining it with Chebychev's inequality yields
\begin{equation}
\sup_{\theta\in\Theta}\bE_\theta(g_n[\bX]-\theta)^2\ge\frac{\sigma^2}{n}\cdot\max_{0\le\alpha
1/2}\alpha(1-2\alpha)=\frac{0.125\sigma^2}{n}.
\end{equation}
In \cite[p.\ 257]{Lehmann83}, it is shown that when $\Theta=\reals$, 
the minimax estimator for this model is the
sample mean, and so, in this case, the correct constant in front of $\sigma^2$ is
actually 1.\\

\noindent
{\bf Example 5.}
Let $X(t)=s(t,\theta)+N(t)$, $t\in[0,T]$, where $N(t)$ is additive white
Gaussian noise (AWGN) with
double-sided spectral density $N_0/2$ and $s(t,\theta)$ is a deterministic signal that
depends on the unknown parameter, $\theta$. 
It is assumed that
the signal energy,
$E=\int_0^Ts^2(t,\theta)\mbox{d}t$, 
does not depend on $\theta$ (which is the case, for example, when $\theta$ is
a delay parameter of a pulse fully contained in the observation interval, or
when $\theta$ is the frequency or the phase of a sinusoidal waveform).
We further assume that $s(t,\theta)$ is at least twice differentiable w.r.t.\ $\theta$, 
and that the energies of the first two derivatives are also independent of
$\theta$. Then, as shown in Appendix A, for small
$|\theta_1-\theta_0|$,
\begin{eqnarray}
\label{approxrho}
\varrho(\theta_0,\theta_1)&\dfn&\frac{1}{E}\int_0^Ts(t,\theta_0)s(t,\theta_1)\mbox{d}t\nonumber\\
&=&1-\frac{(\theta_1-\theta_0)^2}{2E}\int_0^T\left[\frac{\partial
s(t,\theta)}{\partial\theta}\bigg|_{\theta=\theta_0}\right]^2\mbox{d}t+o(|\theta_1-\theta_0|^2)\nonumber\\
&=&1-\frac{(\theta_2-\theta_1)^2\dot{E}}{2E}+o(|\theta_1-\theta_0|^2,
\end{eqnarray}
where $\dot{E}$ is the energy $\dot{s}(t,\theta)=\partial s(t,\theta)/\partial\theta$.

The optimal LRT in deciding between the two hypotheses is based on comparing
between the correlations, $\int_0^TX(t)s(t,\theta_0)\mbox{d}t$ and
$\int_0^TX(t)s(t,\theta_1)\mbox{d}t$. Again, the optimal value of $q$ is
$q^*=1/2$. Thus,
\begin{eqnarray}
P_{\mbox{\tiny e}}\left(\frac{1}{2},\theta_0,\theta_1\right)&=&
Q\left(\sqrt{\frac{E}{N_0}\left[1-\varrho(\theta_0,\theta_1)\right]}\right)\nonumber\\
&\approx&Q\left(\sqrt{\frac{E\dot{E}(\theta_1-\theta_0)^2}{2EN_0}}\right)\nonumber\\
&=&Q\left(\sqrt{\frac{\dot{E}}{2N_0}}\cdot|\theta_1-\theta_0|\right)\nonumber\\
&=&Q\left(\sqrt{\frac{\dot{P}}{2N_0}}\cdot\sqrt{T}|\theta_1-\theta_0|\right),
\end{eqnarray}
where $\dot{P}=\dot{E}/T$ is the power of $\dot{s}(t,\theta)$.
Since we are dealing here with continuous time, then instead of a sequence
$\xi_n$, we use a function, $\xi(T)$, of the observation time, $T$, which
in the case would be $\xi(T)=\frac{1}{\sqrt{T}}$. Let $\theta_0=\theta$ and
$\theta_1=\theta+\frac{2s}{\sqrt{T}}$. Then,
\begin{equation}
P_{\mbox{\tiny e}}^\infty(\theta,s)=Q\left(\sqrt{\frac{2\dot{P}}{N_0}}\cdot |s|\right),
\end{equation}
which, for the MSE case, yields
\begin{eqnarray}
r(\theta)&\ge&\sup_{s\ge 0}\left\{2s^2Q\left(\sqrt{2\dot{P}}{N_0}\cdot
s\right)\right\}\nonumber\\
&=&\frac{N_0}{\dot{P}}\sup_{u\ge 0}\{u^2Q(u)\}\nonumber\\
&=&\frac{0.1657N_0}{\dot{P}}.
\end{eqnarray}
w.r.t.\ $\zeta(T)=1/(1/\sqrt{T})^2=T$, which means that the minimax loss is
lower bounded by $r(\theta)/T\ge 0.1657N_0/\dot{E}$.
This has the same form as the Cram\'er-Rao lower bound (CRLB), except that the multiplicative factor is
0.1657 rather than 0.5. It should be kept in mind, however, that the CRLB is
limited to unbiased estimators.

In \cite[eq.\ (20)]{ZZ69}, the bound is of the same form, but with a multiplicative constant
of 0.16 (for large signal-to-noise ratio), but it should be kept in mind that their notation of the double-sided spectral
density of the noise is $N_0$, rather than $N_0/2$ as here, and so for a fair comparison,
their constant should actually be replaced by $0.16/2=0.08$.\\

The case where $s(t,\theta)$ is not everywhere differentiable w.r.t.\ $\theta$
can be handled in a similar manner, but some caution should be exercised.
For example, consider the model, 
\begin{equation}
X(t)=s(t-\theta)+N(t), 
\end{equation}
where $-\infty < t < \infty$, $\theta\in[0,T]$, $N(t)$ is AWGN as before, and
$s(\cdot)$ is a rectangular pulse
with duration $\Delta$ and amplitude $\sqrt{E/\Delta}$, $E$ being the
signal energy. 
Here, $\varrho(\theta,\theta+\delta)=1-|\delta|/\Delta$, namely, it includes
also a linear term in $|\delta|$, not just the quadratic one. This changes the
asymptotic behavior of the resulting lower bound to $r(\theta)$ turns out to be
$0.1886(N_0\Delta/P)^2$ w.r.t.\ $T^2$ (namely, a minimax
lower bound of $0.1886(N_0\Delta/E)^2$).
It is interesting to compare this bound to the Chapman-Robbins bound for the
same model, which is a local bound of the same form but with a multiplicative constant of
$0.0405$ instead of $0.1886$, and which is limited to unbiased estimators.\\

\noindent
{\bf Example 6.}
Consider an exponential family, 
\begin{equation}
p(\bx|\theta)=\prod_{i=1}^np(x_i|\theta)=\prod_{i=1}^n\frac{e^{\theta
T(x_i)}}{Z(\theta)}=\frac{\exp\left\{\theta\sum_{i=1}^nT(x_i)\right\}}{Z^n(\theta)}.
\end{equation}
where $T(\cdot)$ is a given function and $Z(\theta)$ is a normalization
function given by
\begin{equation}
Z(\theta)=\int_{\reals}e^{\theta T(x)}\mbox{d}x.
\end{equation}
In the binary hypothesis problem, the test statistic is $\sum_{i=1}^nT(X_i)$.
If $q=1/2$, $\xi_n=1/\sqrt{n}$ and $\theta_1-\theta_0=2s/\sqrt{n}$, the LRT
amounts to examining
whether $\sum_{i=1}^n[T(X_i)-\bE_\theta\{T(X_i)\}]$ is larger than
$$\frac{s}{\sqrt{n}}\cdot n\frac{\mbox{d}^2\ln Z(\theta)}{\mbox{d}\theta^2}=
s\sqrt{n}\cdot \frac{\mbox{d}^2\ln Z(\theta)}{\mbox{d}\theta^2}.$$
In this case, the probability of error can be asymptotically assessed using
the central limit theorem (CLT), which after a simple algebraic manipulation,
becomes:
\begin{eqnarray}
P_{\mbox{\tiny e}}^\infty(\theta,s)&=&
Q\left(\frac{s\mbox{d}^2\ln Z(\theta)/\mbox{d}\theta^2}{\sqrt{\mbox{d}^2\ln
Z(\theta)/\mbox{d}\theta^2}}\right)\nonumber\\
&=&Q\left(s\sqrt{\mbox{d}^2\ln Z(\theta)/\mbox{d}\theta^2}
\right)\nonumber\\
&=&Q(s\sqrt{I(\theta)}),
\end{eqnarray}
where $I(\theta)$ is the Fisher information.
Thus, for the MSE,
\begin{equation}
r(\theta)\ge \sup_{s\ge 0} 2s^2Q\left(s\sqrt{I(\theta)}\right)=\frac{0.3314}{I(\theta)}
\end{equation}
w.r.t.\ $\zeta_n=1/(1/\sqrt{n})^2=n$.
In the case of an exponential family with a $d$-dimensional parameter vector,
and $\rho(\varepsilon)=\|\varepsilon\|^2$,
this derivation extends to yield $r(\theta)\ge 0.3314/\lambda_{\min}(\theta)$,
where $\lambda_{\min}(\theta)$ denotes the smallest eigenvalue of the
$d\times d$ Fisher information matrix, $I(\theta)=\nabla^2\ln Z(\theta)$.
We omit the details of this derivation.

\subsection{Extensions of Theorem 1 to Multiple Test Points}

In this subsection, we present several extensions of Theorem 1, from two test
points, $\theta_0$ and $\theta_1$, to a general number, $m$, of test points,
$\theta_0,\theta_1,\ldots,\theta_{m-1}$, with corresponding priors,
$q_0,q_1,\ldots,q_{m-1}$, which are non-negative reals that sum to unity.
Some of these bounds may induce corresponding extensions of Corollary 1
by letting $\theta_0,\theta_1,\ldots,\theta_{m-1}$ approach each other as $n$
grows, so that the corresponding error probabilities would converge to
positive constants.
We have not examined these extensions numerically. For $m$ test points, the
number of degrees of freedom to optimize in order to calculate the tightest
possible bound is $m(d+1)-1$, which consists of the $m$ $d$-dimensional
parameters, $\theta_0,\ldots,\theta_{m-1}$ plus the $m-1$ degrees of freedom
associated with $q_0,\ldots,q_{m-2}$.

The first step, that is common to all extended versions to be presented in
this subsection, is to bound $R_n(\Theta)$ from below by
\begin{equation}
R_n(\Theta)\ge
\inf_{g_n}\sup_{\theta_0,\ldots,\theta_{m-1}}\sup_{q_0,\ldots,q_{m-2}}
\sum_{0=1}^{m-1}q_i\bE_{\theta_i}\{\rho(g_n[\bX]-\theta_i)\}. 
\end{equation}
\\

\noindent
{\it A.~Bounds Based on Pairwise Terms}\\

The first approach is to
further manipulate the weighted sum on the basis of
pairwise terms similarly as in Theorem 1. 

One such resulting inequality is
the following:
\begin{eqnarray}
\sum_{i=0}^{m-1}q_i\bE_{\theta_i}\{\rho(\hat{\theta}-\theta_i)\}&=&
\frac{1}{2}\sum_{i=0}^{m-1}[q_i\bE_{\theta_i}\{\rho(\hat{\theta}-\theta_i)\}+q_{i\oplus 1}
\bE_{\theta_{\oplus 1}i}\{\rho(\hat{\theta}-\theta_{i\oplus 1})\}]\nonumber\\
&\ge&\sum_{i=0}^{m-1} \rho\left(\frac{\theta_{i\oplus
1}-\theta_i}{2}\right)(q_i+q_{i\oplus 1})P_{\mbox{\tiny
e}}\left(\frac{q_i}{q_i+q_{i\oplus 1}},\theta_i,\theta_{i\oplus 1}\right),
\end{eqnarray}
where $\oplus$ denotes addition modulo $m$ (so that $(m-1)\oplus 1=0$),
and where the last inequality is proved similarly as in the proof of Theorem 1.
Another option is to apply the same idea to all possible pairs of test points,
i.e.,
\begin{eqnarray}
\sum_{i=0}^{m-1}q_i\bE_{\theta_i}\{\rho(\hat{\theta}-\theta_i)\}&=&
\frac{1}{2(m-1)}\sum_{i\ne
j}[q_i\bE_{\theta_i}\{\rho(\hat{\theta}-\theta_i)\}+q_j\bE_{\theta_j}\{\rho(\hat{\theta}-\theta_j)\}]\nonumber\\
&\ge&\frac{1}{m-1}\sum_{i\ne j}\rho\left(\frac{\theta_j
-\theta_i}{2}\right)(q_i+q_j)P_{\mbox{\tiny
e}}\left(\frac{q_i}{q_i+q_j},\theta_i,\theta_j\right).
\end{eqnarray}
This first version is more suitable to a one-dimensional parameter, where the values of
$\{\theta_i\}$ form a certain grid along a line (see, e.g., \cite{me}).
The second version is
more natural in the case of a vector parameter, where the terms corresponding
to the various pairs
$(\theta_i,\theta_j)$ explore various directions in the parameter space.
Note that for a given estimator, $g_n[\cdot]$, the weight vector $\bq=(q_0,q_1,\ldots,q_{m-1})$,
that maximizes 
$\sum_{i=0}^{m-1}q_i\bE_{\theta_i}\{\rho(\hat{\theta}-\theta_i)\}$, puts all
its mass on the value (or values) of $i\in\{0,1,\ldots,m-1\}$ with the largest expected loss,
$\bE_{\theta_i}\{\rho(\hat{\theta}-\theta_i)\}$.
This is due to the linearity of the weighted sum in $\bq$. However,
this is not necessarily true when it comes to the lower bounds, as they are no
longer linear functions of $\bq$. In fact, they are concave functions of
$\bq$.\\

\noindent
{\it B.~Bounds Based on Unitary Transformations and the List Error
Probability}\\

Let $\theta$ be a parameter vector of
dimension $d$ and $\Theta\subseteq\reals^d$.
Let $\rho(\varepsilon)$ be a convex loss function 
that depends on the $d$-dimensional error vector $\varepsilon$ only via
its norm, $\|\varepsilon\|^2$ (that is, $\rho$ has radial
symmetry). Let $T_0,T_1,T_2,\ldots,T_{m-1}$ be a set of $m-1$ unitary
transformation
matrices with the properties: (i) $\theta\in\Theta$ if and only if
$T_i\theta\in\Theta$ for all $i=0,1,\ldots,m-1$, and (ii) $T_0+T_1+\cdots+T_{m-1}=0$.
For example, if $d=2$
and $m=3$, take $T_i$ to be matrices of rotation by $2\pi i/3$ ($i=0,1,2$).
More generally, we may
allow also non-linear transformations, $\{T_i(\cdot), i=0,1,\ldots,m-1\}$ that
all map $\Theta$ onto itself, all preserve norms,
and that $\sum_{i=0}^{m-1}T_i(\theta)=0$ for all $\theta\in\Theta$.\\

\noindent
{\bf Theorem 2.}
Let $\theta_0,\theta_1,\ldots,\theta_{m-1}$, $q_0,q_1,\ldots,q_{m-1}$ be
given, and let
$T_0,T_1,\ldots T_{m-1}$ be as described in the above paragraph. Then,
for a convex loss function $\rho(\varepsilon)$ that depends on $\varepsilon$
only via $\|\varepsilon\|$, we have:
\begin{equation}
R_n(\Theta)\ge
m\cdot\rho\left(\frac{1}{m}\sum_{i=0}^{m-1}T_i\theta_i\right)\cdot
\int_{\reals^n}
\min\{q_0\cdot
p(\bx|\theta_0),\ldots,q_{m-1}p(\bx|\theta_{m-1})\}\mbox{d}\bx.
\end{equation}

The integral on the right-hand side can be interpreted as the probability of
list-error with a list size of $m-1$. In other words, imagine a multiple hypothesis
testing problem where the observer, upon observing $\bx$, constructs a list of the $m-1$ most
likely hypotheses in descending order of $q_ip(\bx|\theta_i)$. Then, the above
integral can be identified as the list-error probability, namely, the
probability that the correct index $i$ is not in the list.\\

\noindent
{\it Proof of Theorem 2.}
The proof is a direct extension of the proof of Theorem 1: 
\begin{eqnarray}
\sup_{\theta\in\Theta}\bE_\theta\{\rho(g_n[\bX]-\theta)\}
&\ge&
\sum_{i=0}^{m-1}q_i\bE_{\theta_i}\{\rho(g_n[\bX]-\theta_i)\}\nonumber\\
&=&\int_{\reals^n}
\bigg[\sum_{i=0}^{m-1}q_i\cdot
p(\bx|\theta_i)\rho(g_n[\bx]-\theta_i)\bigg]\mbox{d}\bx\nonumber\\
&\ge&
m\cdot\int_{\reals^n}
\min\{q_0\cdot p(\bx|\theta_0),\ldots,q_{m-1}p(\bx|\theta_{m-1})\}\times\nonumber\\
& &\bigg[\frac{1}{m}\sum_{i=0}^{m-1}\rho(\theta_i-g_n[\bX])\bigg]
\mbox{d}\bx\nonumber\\
&=& m\cdot
\int_{\reals^n}
\min\{q_0\cdot p(\bx|\theta_0),\ldots,q_{m-1}p(\bx|\theta_{m-1})\}\times\nonumber\\
& &\bigg[\frac{1}{m}\sum_{i=0}^{m-1}\rho(T_i(\theta_i-g_n[\bX]))\bigg]
\mbox{d}\bx\nonumber\\
&\ge& m\cdot
\int_{\reals^n}
\min\{q_0\cdot p(\bx|\theta_0),\ldots,q_{m-1}p(\bx|\theta_{m-1})\}\times\nonumber\\
& &\rho\left(\frac{1}{m}\sum_{i=0}^{m-1}T_i(\theta_i-g_n[\bX])\right)
\mbox{d}\bx\nonumber\\
&=& m\cdot
\int_{\reals^n}
\min\{q_0\cdot p(\bx|\theta_0),\ldots,q_{m-1}p(\bx|\theta_{m-1})\}\times\nonumber\\
& &\rho\left(\frac{1}{m}\sum_{i=0}^{m-1}T_i\theta_i\right)
\mbox{d}\bx\nonumber\\
&=&m\rho\left(\frac{1}{m}\sum_{i=0}^{m-1}T_i\theta_i\right)\cdot
\int_{\reals^n}
\min\{q_0\cdot
p(\bx|\theta_0),\ldots,q_{m-1}p(\bx|\theta_{m-1})\}\mbox{d}\bx.\nonumber
\end{eqnarray}
This completes the proof of Theorem 2.\\

Theorem 1 is a special case where $m=2$, $T_0=I$ and $T_1=-I$, where $I$
is the $d\times d$ identity matrix. The integral associated with the lower
bound of Theorem 2
might not be trivial to evaluate in general for $m\ge 3$. 
However, there are some choices of the auxiliary parameters that may
facilitate calculations. One such choice is follows.
For some positive integer $k\le m$, take
$\theta_0=\theta_1=\ldots=\theta_{k-1}\dfn \vartheta_0$,
for some $\vartheta_0\in\Theta$,
$q_0=q_1=\ldots=q_{k-1}\dfn Q/k$ for some $Q\in(0,1)$,
$\theta_k=\theta_{k+1}=\ldots=\theta_{m-1}\dfn\vartheta_1$, for some
$\vartheta_1\in\Theta$, and finally,
$q_k=q_{k+1}=\ldots=q_{m-1}=(1-Q)/(m-k)$.
The integrand then becomes the minimum between two
functions only, as before. Denoting $\alpha=k/m$, the bound then becomes
\begin{eqnarray}
R_n(\Theta)
&\ge&m\rho\left(\frac{1}{m}\sum_{i=0}^{k-1}T_i(\vartheta_0-\vartheta_1)\right)\cdot\int_{\reals^n}
\min\left\{\frac{Q}{k}\cdot p(\bx|\vartheta_0),\frac{1-Q}{m-k}\cdot p(\bx|\vartheta_1)\right\}\mbox{d}\bx\nonumber\\
&=&\rho\left(\frac{1}{m}\sum_{i=0}^{k-1}T_i(\vartheta_0-\vartheta_1)\right)\cdot\int_{\reals^n}
\min\left\{\frac{Q}{\alpha}\cdot p(\bx|\vartheta_0),\frac{1-Q}{1-\alpha}\cdot p(\bx|\vartheta_1)\right\}\mbox{d}\bx\nonumber\\
&=&\rho\left(\frac{1}{m}\sum_{i=0}^{k-1}T_i(\vartheta_0-\vartheta_1)\right)
\cdot\left(\frac{Q}{\alpha}+\frac{1-Q}{1-\alpha}\right)\cdot P_{\mbox{\tiny
e}}\left(\frac{(1-\alpha)Q}{(1-\alpha)Q+\alpha(1-Q)},\vartheta_0,\vartheta_1\right).
\end{eqnarray}
Redefining
\begin{equation}
q=\frac{(1-\alpha)Q}{(1-\alpha)Q+\alpha(1-Q)},
\end{equation}
we have 
\begin{equation}
\frac{Q}{\alpha}+\frac{1-Q}{1-\alpha}=\frac{1}{1-\alpha-q+2\alpha q},
\end{equation}
and the following corollary to Theorem 2 is obtained.\\

\noindent
{\bf Corollary 2.} 
Let the conditions of Theorem 2 be satisfiled. Then,
\begin{equation}
R_n(\Theta)\ge\sup_{\vartheta_0,\vartheta_1,\alpha,q}
\rho\left(\frac{1}{m}\sum_{i=0}^{k-1}T_i(\vartheta_0-\vartheta_1)\right)\cdot\frac{P_{\mbox{\tiny
e}}(q,\vartheta_0,\vartheta_1)}{1-\alpha-q+2\alpha q}.
\end{equation}
\\

Note that if $m$ is even, $\alpha=\frac{1}{2}$ and $T_0=I$, then we are actually back
to the bound of $m=2$, and so, 
the optimal bound for even $m > 2$
cannot be worse than our bound of $m=2$. We do not have, however, precisely the same
argument for odd $m$, but for large $m$ it becomes immaterial if $m$ is even
or odd. In its general form, the bound of Theorem 2 is a heavy optimization problem, as we
have the freedom to optimize $\theta_0,\ldots,\theta_{m-1}$,
$T_0,\ldots,T_{m-1}$ (under the constraints that they are all unitary and sum
to zero), and $q_0,\ldots,q_{m-1}$ (under the constraints that they are all
non-negative and sum to unity). 

Another relatively convenient choice is to take
$\theta_i=T_i^{-1}\theta_0$, $i=1,\ldots,m-1$, to obtain another corollary to
Theorem 2:\\

\noindent
{\bf Corollary 3.} Let the conditions of Theorem 2 be satisfied. Then,
\begin{equation}
R_n(\Theta)\ge\sup_{\theta_0,T_0,\ldots,T_{m-1},q_0,\ldots,q_{m-1}}
m\rho(\theta_0)\cdot\int_{\reals^n}\min\{q_0p(\bx|\theta_0),q_1p(\bx|T_1^{-1}\theta_0),\ldots,
q_{m-1}p(\bx|T_{m-1}^{-1}\theta_0)\}\mbox{d}\bx.
\end{equation}

\noindent
{\bf Example 7.}
To demonstrate a calculation of the extended lower bound for $m=3$,
consider the following model. We are observing a noisy signal,
\begin{equation}
Z_i=\vartheta\phi_i+(\vartheta+\zeta)\psi_i+N_i,~~~~~i=1,2,\ldots,n,
\end{equation}
where $\vartheta$ is the desired parameter to be estimated, $\zeta$ is a nuisance
parameter, taking values within an interval $[-\delta,\delta]$
for some $\delta>0$, 
$\{N_i\}$ are i.i.d.\ Gaussian random variables with zero mean and variance
$\sigma^2$, and 
$\phi_i$ and $\psi_i$ are two given orthogonal waveforms with
$\sum_{i=1}^n\phi_i^2=
\sum_{i=1}^n\psi_i^2=n$. Suppose we are interested in estimating $\vartheta$
based on the sufficient statistics
$X=\frac{1}{n}\sum_{i=1}^nZ_i\phi_i$
and $Y=\frac{1}{n}\sum_{i=1}^nZ_i\psi_i$, which are
jointly Gaussian random variables with mean vector $(\vartheta,\vartheta+\zeta)$ and covariance matrix
$\frac{\sigma^2}{n}\cdot I$, $I$ being the $2\times 2$ identity matrix.
We denote realizations of $(X,Y)$ by $(x,y)$. 
Let us also denote $\theta=(\vartheta,\vartheta+\zeta)$. Since we are interested
only in estimating $\vartheta$, our loss function will depend only on the
estimation error of the first component of $\theta$, which is $\vartheta$.
Consider the choice $m=3$
and let $T_i$ be counter-clockwise rotation
transformations by $2\pi i/3$, $i=0,1,2$.
For a given $\Delta\in(0,\delta]$, let us select $\theta_0=(-\Delta,0)$,
$\theta_1=T_1^{-1}\theta_0=(\Delta/2,\Delta\sqrt{3}/2)$ and
$\theta_2=T_2^{-1}\theta_0=(\Delta/2,-\Delta\sqrt{3}/2)$. Finally, let
$q_0=q_1=q_2=\frac{1}{3}$. In order to calculated the integral
$$I=\int_{\reals^2}\min\left\{\frac{1}{3}p(x,y|\theta_0),
\frac{1}{3}p(x,y|\theta_1),\frac{1}{3}p(x,y|\theta_2)\right\}\mbox{d}x\mbox{d}y$$
the plane $\reals^2$ can be partitioned into three slices over which the integrals
contributed are equal. In each such region, the smallest
$p(x,y|\theta_i)$ is integrated. In other words, every
$p(x,y|\theta_i)$ in its turn is integrated over the region whose
Euclidean distance to $\theta_i$ is larger than the distances to the other two values
of $\theta$. For $\theta_0=(-\Delta,0)$, this is the region $\{(x,y):~x\ge 0,~|y|\le x\sqrt{3}\}$.
The factor of $\frac{1}{3}$ cancels with the three identical contributions
from $\theta_0$, $\theta_1$ and $\theta_2$ due to the symmetry. Therefore,
\begin{eqnarray}
I&=&\int_0^\infty\mbox{d}x\int_{-x\sqrt{3}}^{x\sqrt{3}}\mbox{d}y\cdot
p(x,y|\theta_0)\nonumber\\
&=&\int_0^\infty\mbox{d}x\int_{-x\sqrt{3}}^{x\sqrt{3}}\cdot\frac{n}{2\pi\sigma^2}
\exp\left\{-\frac{n(x+\Delta)^2+ny^2}{2\sigma^2}\right\}\mbox{d}y\nonumber\\
&=&\sqrt{\frac{n}{2\pi\sigma^2}}\cdot\int_0^\infty
\exp\left\{-\frac{n(x+\Delta)^2}{2\sigma^2}\right\}\left[1-2Q\left(\frac{x\sqrt{3n}}{\sigma}\right)\right]\mbox{d}x.
\end{eqnarray}

Our next mathematical manipulations in this example are in the spirit of the passage from Theorem 1 to
Corollary 1, that is, selecting the test points increasingly close to each
other as functions of $n$, so that
the probability of list-error would tend to a positive constant. To this end,
we change the integration variable $x$ to $u=x\sqrt{n}/\sigma$ and select
$\Delta=s\sigma/\sqrt{n}$ for some $s\ge 0$ to be optimized later. Then,
\begin{eqnarray}
I&=&\sqrt{\frac{n}{2\pi\sigma^2}}\cdot\int_0^\infty
\exp\left\{-\frac{n(u\sigma/\sqrt{n}+s\sigma/\sqrt{n})^2}{2\sigma^2}\right\}\left[1-2Q(u\sqrt{3})\right]
\mbox{d}\left(\frac{u\sigma}{\sqrt{n}}\right)\nonumber\\
&=&\frac{1}{\sqrt{2\pi}}\cdot\int_0^\infty
e^{-(u+s)^2/2}\cdot[1-2Q(u\sqrt{3})]
\mbox{d}u.
\end{eqnarray}
The MSE bound then becomes
\begin{eqnarray}
R_n(\theta,g_n)&\ge&
\sup_{s\ge 0}\left\{3\cdot\left(\frac{s\sigma}{\sqrt{n}}\right)^2\cdot
\frac{1}{\sqrt{2\pi}}\cdot\int_0^\infty
e^{-(u+s)^2/2}\cdot[1-2Q(u\sqrt{3})]
\mbox{d}u\nonumber\right\}\\
&=&\frac{\sigma^2}{n}\cdot\sup_{s\ge 0}\left\{\frac{3s^2}{\sqrt{2\pi}}\cdot
\int_0^\infty
e^{-(u+s)^2/2}\cdot[1-2Q(u\sqrt{3})]
\mbox{d}u\right\}\nonumber\\
&=&\frac{0.2514\sigma^2}{n}.
\end{eqnarray}
This bound is not as tight as the corresponding bound of $m=2$, which results
in $0.3314\sigma^2/n$, but it should be kept in mind that here, we have not
attempted to optimize the choices of $\theta_0$, $\theta_1$, $\theta_2$, $T_0$, $T_1$,
$T_2$, $q_0$, and $q_1$. Instead, we have chosen these parameter values from
considerations of computational convenience, just to demonstrate the
calculation. This concludes Example 7.\\

\section{Bounds Based on the Minimum Expected Loss Over Some Test Points}
\label{moments}

\subsection{Two Test Points}

The following generic, yet conceptually very simple, lower bound assumes neither symmetry, nor convexity 
of the loss function $\rho(\cdot)$. For a given
$(\bx,q,\theta_0,\theta_1)\in\reals^n\times[0,1]\times\Theta^2$,
let us define
\begin{equation}
\psi(\bx,q,\theta_0,\theta_1)=\min_{\vartheta}\{qp(\bx|\theta_0)\rho(\vartheta-\theta_0)+
(1-q)p(\bx|\theta_1)\rho(\vartheta-\theta_1)\}.
\end{equation}
Then,
\begin{eqnarray}
\label{generic}
R_n(\Theta)&\ge&
\sup_{\theta_0,\theta_1,q}q\bE_{\theta_0}\{\rho(g_n[\bX]-\theta_0)\}+
(1-q)\bE_{\theta_1}\{\rho(g_n[\bX]-\theta_1)\}\nonumber\\
&=&\sup_{\theta_0,\theta_1,q}\int_{\reals^n}[qp(\bx|\theta_0)\rho(g_n[\bx]-\theta_0)+
(1-q)p(\bx|\theta_1)\rho(g_n[\bx]-\theta_1)]\mbox{d}\bx\nonumber\\
&\ge&\sup_{\theta_0,\theta_1,q}\int_{\reals^n}\psi(\bx,q,\theta_0,\theta_1)\mbox{d}\bx.
\end{eqnarray}
If we further assume symmetry of $\rho$, then it is easy to see that the minimizer
$\vartheta^*$, that achieves $\psi(\bx,q,\theta_0,\theta_1)$, is always within
the interval $[\theta_0,\theta_1]$. This is because the objective increases
monotonically as we move away from the interval $[\theta_0,\theta_1]$
in either direction.
Of course, this simple idea can
easily be extended to apply to weighted sums of more than two points, in
principle, but it
would become more complicated -- see the next subsection for three such points.

If $\rho$ is concave, then the minimizing $\vartheta$ is either $\theta_0$ or
$\theta_1$, depending on the smaller between $qp(\bx|\theta_0)$ and
$(1-q)p(\bx|\theta_1)$ and the bound becomes
\begin{equation}
R_n(\Theta)\ge \sup_{\theta_0,\theta_1,q}\rho(\theta_1-\theta_0)P_{\mbox{\tiny e}}(q,\theta_0,\theta_1).
\end{equation}
Minimality at the edge-points may happen also for some loss functions that are
not concave, like the loss function $\rho(u)=1\{|u|\ge\Delta\}$.

The generic lower bound (\ref{generic}) is more general than our first bound in the sense that it does not
require convexity or symmetry of $\rho$, but the down side is that the resulting
expressions are harder
to deal with directly, as will be seen shortly. For loss
functions other than the MSE or general moments of the estimation error, it may not be a trivial task even to derive a
closed form expression of $\psi(\bx,q,\theta_0,\theta_1)$ (i.e., to carry out
the minimization associated with the definition of $\psi$).

For the case of the MSE, $\rho(\varepsilon)=\varepsilon^2$, 
the calculation of $\psi(\bx,q,\theta_0,\theta_1)$ is straightforward, and it
readily yields
\begin{equation}
\psi(\bx,q,\theta_0,\theta_1)=(\theta_1-\theta_0)^2\cdot\frac{qp(\bx|\theta_0)\cdot(1-q)p(\bx|\theta_1)}{qp(\bx|\theta_0)+
(1-q)p(\bx|\theta_1)}.
\end{equation}
However, it may not be convenient to integrate this function of $\bx$ due to
the summation at the denominator. One way to alleviate this difficulty is to
observe that
\begin{eqnarray}
\label{37}
\psi(\bx,q,\theta_0,\theta_1)&\ge&
(\theta_1-\theta_0)^2\cdot\frac{qp(\bx|\theta_0)\cdot(1-q)p(\bx|\theta_1)}{2\max\{qp(\bx|\theta_0),
(1-q)p(\bx|\theta_1)\}}\nonumber\\
&=&\frac{1}{2}\cdot
(\theta_1-\theta_0)^2\cdot\min\{qp(\bx|\theta_0),(1-q)p(\bx|\theta_1)\},
\end{eqnarray}
which after integration yields again,
\begin{equation}
R_n(\Theta)\ge\sup_{\theta_0,\theta_1,q}\left\{\frac{1}{2}(\theta_1-\theta_0)^2\cdot
P_{\mbox{\tiny e}}(q,\theta_0,\theta_1)\right\},
\end{equation}
exactly as in Theorem 1 in the special case of the MSE. This indicates that
the bound (\ref{generic}) is at least as tight as the bound of Theorem 1 for
the MSE.

It turns out, however, that we can do better than bounding the denominator, $qp(\bx|\theta_0)+
(1-q)p(\bx|\theta_1)$, by $2\cdot\max\{qp(\bx|\theta_0),
(1-q)p(\bx|\theta_1)\}$ for the purpose of obtaining a more convenient
integrand. Specifically, consider the identity stated following lemma, whose
proof appears in Appendix B.

\noindent
{\bf Lemma 1.}
Let $k$ be a positive integer and let $a_1,\ldots,a_k$ be positive reals.
Then,
\begin{equation}
\sum_{i=1}^ka_i
=\inf_{(r_1,\ldots,r_k)\in\calS}\max\left\{\frac{a_1}{r_1},\ldots,\frac{a_k}{r_k}\right\},
\end{equation}
where $\calS$ is the interior of the $k$-dimensional simplex, namely, the set
of all vectors
$(r_1,\ldots,r_k)$ with strictly positive components that sum to unity.\\

Applying Lemma 1 with the assignments $k=2$, $a_1=qp(\bx|\theta_0)$ and
$a_2=(1-q)p(\bx|\theta_1)$, we have
\begin{equation}
qp(\bx|\theta_0+(1-q)p(\bx|\theta_1)
=\inf_{r\in(0,1)}
\max\left\{\frac{qp(\bx|\theta_1)}{r},\frac{(1-q)p(\bx|\theta_2)}{1-r}\right\}.
\end{equation}
Thus,
\begin{eqnarray}
\psi(\bx,q,\theta_0,\theta_1)&=&\frac{qp(\bx|\theta_1)\cdot(1-q)p(\bx|\theta_2)}{\inf_{r\in(0,1)}
\max\left\{qp(\bx|\theta_0)/r,(1-q)p(\bx|\theta_1)/(1-r)\right\}}\nonumber\\
&=&\sup_{r\in(0,1)}\frac{qp(\bx|\theta_0)\cdot(1-q)p(\bx|\theta_1)}{
\max\left\{qp(\bx|\theta_0)/r,(1-q)p(\bx|\theta_1)/(1-r)\right\}}\nonumber\\
&=&\sup_{r\in(0,1)}\left\{\begin{array}{ll}
r(1-q)p(\bx|\theta_1) & r \le r^*\\
(1-r)qp(\bx|\theta_0) & r \ge r^*\end{array}\right.\nonumber\\
&=&\sup_{r\in(0,1)}\min\{r(1-q)p(\bx|\theta_1),(1-r)qp(\bx|\theta_0)\},
\end{eqnarray}
where 
\begin{equation}
r^*=\frac{qp(\bx|\theta_0)}{qp(\bx|\theta_0)+(1-q)p(\bx|\theta_1)}.
\end{equation}
Thus, the bound becomes
\begin{eqnarray}
R_n(\Theta)&\ge&
\sup_{\{(\theta_0,\theta_1,q)\in\Theta^2\times(0,1)\}}
(\theta_1-\theta_0)^2\times\nonumber\\
& &\int_{\reals^n}\sup_{r\in(0,1)}\min\{r(1-q)p(\bx|\theta_1),(1-r)qp(\bx|\theta_0)\}\mbox{d}\bx\nonumber\\
&\ge&\sup_{\{(\theta_0,\theta_1,q)\in\Theta^2\times(0,1)^2\}}\sup_{r\in(0,1)}
(\theta_1-\theta_0)^2\times\nonumber\\
& &\int_{\reals^n}\min\{r(1-q)p(\bx|\theta_1),(1-r)qp(\bx|\theta_0)\}\mbox{d}\bx\nonumber\\
&=&\sup_{\{(\theta_0,\theta_1,q,r)\in\Theta^2\times(0,1)^2\}}
(\theta_1-\theta_0)^2\cdot(q+r-2qr)\times\nonumber\\
& &\int_{\reals^n}\min\left\{\frac{r(1-q)}{q+r-2qr}\cdot p(\bx|\theta_1),
\frac{(1-r)q}{q+r-2qr}\cdot p(\bx|\theta_0)\right\}\mbox{d}\bx\nonumber\\
&=&\sup_{\{(\theta_1,\theta_2,q,r)\in\Theta^2\times(0,1)^2\}}
(\theta_1-\theta_0)^2\cdot(q+r-2qr)\cdot
P_{\mbox{\tiny
e}}\left(\frac{(1-r)q}{q+r-2qr},\theta_0,\theta_1\right).\nonumber
\end{eqnarray}
The bound of Theorem 1 for the MSE is obtained as a
special case of
$r=1/2$. Therefore, after the optimization over the additional degree of
freedom, $r$, the resulting bound cannot be worse than the MSE bound of Theorem 1.
In fact, it may strictly improve as we will demonstrate shortly.
The choice $r=q$ gives a prior of $1/2$ in the error probability
factor, and then the maximum of the external factor,
$q+r-2rq=2q(1-q)$, is maximized by $q=1/2$.\\

\noindent
{\bf Example 2'.}
To demonstrate the new bound for the MSE, let us revisit Example 2 and see how it
improves the multiplicative constant. In that example,
\begin{equation}
P_{\mbox{\tiny e}}\left(\frac{(1-r)q}{q+r-2qr},\theta_0,\theta_1\right)=
\min\left\{\frac{(1-r)q}{q+r-2qr},\frac{r(1-q)}{q+r-2qr}\cdot\left(\frac{\theta_0}{\theta_1}\right)^n\right\}.
\end{equation}
Let us denote $\alpha=(\theta_0/\theta_1)^n$ and recall that $\alpha\in(0,1)$,
provided that we select $\theta_1>\theta_0$. Then,
\begin{equation}
(q+r-2qr)\cdot
P_{\mbox{\tiny e}}\left(\frac{(1-r)q}{q+r-2qr},\theta_0,\theta_1\right)=
\min\{(1-r)q,r(1-q)\alpha\}.
\end{equation}
The maximum w.r.t.\ $q$ is attained when $(1-r)q=r(1-q)\alpha$, namely,
for $q=q^*\dfn\alpha r/[1-(1-\alpha)r]$, which yields
\begin{equation}
\max_q(q+r-2qr)\cdot
P_{\mbox{\tiny e}}\left(\frac{(1-r)q}{q+r-2qr},\theta_0,\theta_1\right)=\frac{\alpha
r(1-r)}{1-(1-\alpha)r}.
\end{equation}
Let us denote $\beta\dfn 1-\alpha\in(0,1)$. The maximum of $r(1-r)/(1-\beta
r)$ is attained for
\begin{equation}
r=r^*\dfn\frac{1-\sqrt{1-\beta}}{\beta}=\frac{1-\sqrt{\alpha}}{1-\alpha}=\frac{1}{1+\sqrt{\alpha}},
\end{equation}
which yields
\begin{eqnarray}
\max_{q,r}(q+r-2qr)\cdot
P_{\mbox{\tiny
e}}\left(\frac{(1-r)q}{q+r-2qr},\theta_0,\theta_1\right)&=&
\sup_{r\in(0,1)}\frac{\alpha
r(1-r)}{1-(1-\alpha)r}\nonumber\\
&=&\alpha\cdot\left(\frac{1-\sqrt{\alpha}}{1-\alpha}\right)^2\nonumber\\
&=&\frac{\alpha}{(1+\sqrt{\alpha})^2}.
\end{eqnarray}
To obtain a local bound in the spirit of Corollary 1,
take $\theta_0=\theta$, $\theta_1=\theta(1+\frac{s}{n\theta})$ which yields
$\alpha=e^{-s/\theta}$ in the limit of
large $n$, and so,
\begin{equation}
r(\theta)\ge
\sup_{s\ge 0}
\frac{s^2e^{-s/\theta}}{(1+e^{-s/[2\theta]})^2}=
\theta^2\cdot\sup_{u\ge 0}u^2\cdot
\frac{e^{-u}}{(1+e^{-u/2})^2}=
0.3102\theta^2,
\end{equation}
w.r.t.\ $\zeta_n=n^2$,
which improves on our earlier bound in Example 2, $r(\theta)\ge 0.2414\theta^2$.
This concludes Example 2'.\\

More generally, for general moments of the estimation error,
a similar derivation yields the following:\\

\noindent
{\bf Theorem 3.}
For $\rho(\varepsilon)=|\varepsilon|^{t}$, $t\ge 1$ (not necessarily an
integer),
\begin{equation}
R_n(\Theta)\ge
\sup_{\theta_0,\theta_1,q,r}
|\theta_1-\theta_0|^t[(1-r)^{t-1}q+r^{t-1}(1-q)]\cdot P_{\mbox{\tiny
e}}\left(\frac{(1-r)^{t-1}q}{(1-r)^{t-1}q+r^{t-1}(1-q)},\theta_0,\theta_1\right).
\end{equation}

Applying the local version of Theorem 3 to Example 2',
we get:
\begin{equation}
r(\theta)\ge\theta^t\cdot
\sup_{s>0}\frac{s^te^{-s}}{(1+e^{-s/t})^t}
\end{equation}
w.r.t.\ $\zeta_n=n^t$.
Changing the optimization variable from $s$ to $\sigma=s/t$, we end up with
\begin{eqnarray}
r(\theta)\ge (\theta t)^t\cdot\left[\sup_{\sigma>0}\frac{\sigma}{e^{\sigma}+1}\right]^t
&=&(0.2785t\theta)^t.
\end{eqnarray}
The factor of $(0.2785t)^t$ should be compared with
$(1/2e)^t=0.1839^t$ of Example 2, pertaining to the choice $q=r=1/2$.
The gap increases exponentially with $t$.
For the maximum likelihood estimator pertaining to Example 2, which is $g_n[\bX]=\max_iX_i$, it is easy to show that
whenever $t$ is integer,
\begin{equation}
\bE_\theta\{|\hat{\theta}-\theta|^t\}=\frac{n!t!\theta^t}{(n+t)!}=\frac{t!\theta^t}{(n+1)(n+2)\cdot\cdot\cdot(n+t)},
\end{equation}
and so, the asymptotic gap is between $t!$ and $(0.2785t)^t$.
Considering the Stirling approximation the ratio between the upper bound and
the lower bound is about $\sqrt{2\pi t}\cdot(1.3211)^t$.

\subsection{Three Test Points}

The idea behind the bounds of Subsection 4.1 can be conceptually extended to be based
on more than two test points, but the resulting
expressions become cumbersome very quickly as the number of test points grows. For three points, 
however, this is still manageable and can provide improved bounds. Let us select the three test points to be
$\theta_0-\Delta$,
$\theta_0$ and $\theta_0+\Delta$ for some $\theta_0$ and $\Delta$, and let us
assign weights $q$, $r$, and $w=
1-q-r$. Consider the bound
\begin{eqnarray}
R_n(\Theta)&\ge&
q\bE_{\theta_0-\Delta}\{\rho(g_n[\bX]-\theta_0+\Delta)\}+
r\bE_{\theta_0}\{\rho(g_n[\bX]-\theta_0)\}+\nonumber\\
& &w\bE_{\theta_0+\Delta}\{\rho(g_n[\bX]-\theta_0-\Delta)\}\nonumber\\
&=&\int_{\reals^n}\bigg[q\cdot
p(\bx|\theta_0-\Delta)\rho(g_n[\bx]-\theta_0+\Delta)+
r\cdot p(\bx|\theta_0)\rho(g_n[\bx]-\theta_0)+\nonumber\\
& &w\cdot p(\bx|\theta_0+\Delta)\rho(g_n[\bx]-\theta_0-\Delta)\bigg]\mbox{d}\bx\nonumber\\
&\ge&\int_{\reals^n}\Psi(\bx,\theta_0,\Delta,q,r)\mbox{d}\bx,
\end{eqnarray}
where
\begin{eqnarray}
\Psi(\bx,\theta_0,\Delta,q,r)&=&\min_{\vartheta}\bigg\{
q\cdot p(\bx|\theta_0-\Delta)\rho(\vartheta-\theta_0+\Delta)+
r\cdot p(\bx|\theta_0)\rho(\vartheta-\theta_0)+\nonumber\\
& &w\cdot p(\bx|\theta_0+\Delta)\rho(\vartheta-\theta_0-\Delta)\bigg\}.
\end{eqnarray}
Considering the case of the MSE,
\begin{eqnarray}
\Psi(\bx,\theta_0,\Delta,q,r)&=&\min_{\vartheta}\bigg\{
q\cdot p(\bx|\theta_0-\Delta)(\vartheta-\theta_0+\Delta)^2+
r\cdot p(\bx|\theta_0)(\vartheta-\theta_0)^2+\nonumber\\
& &w\cdot p(\bx|\theta_0+\Delta)(\vartheta-\theta_0-\Delta)^2\bigg\} 
\end{eqnarray}
can be found in closed-form. Denoting temporarily
$a=qp(\bx|\theta_0-\Delta)$,
$b=rp(\bx|\theta_0)$, and 
$c=wp(\bx|\theta_0+\Delta)$, 
$\Psi(\bx,\theta_0,\Delta,q,r)$ is attained by
\begin{equation}
\vartheta=\vartheta^*=\frac{a(\theta_0-\Delta)+b\theta_0+c(\theta_0+\Delta)}{a+b+c}=\theta_0+\frac{(c-a)\Delta}{a+b+c}.
\end{equation}
On substituting $\vartheta^*$ into the sum of squares we end up with
\begin{eqnarray}
\Psi(\bx,\theta_0,\Delta,q,r)&=&\frac{[a(b+2c)^2+b(c-a)^2+c(2a+b)^2]\Delta^2}{(a+b+c)^2}\nonumber\\
&=&\frac{(ab+bc+4ac)\Delta^2}{a+b+c}\nonumber\\
&=&\frac{[qrp(\bx|\theta_0-\Delta)p(\bx|\theta_0)+rwp(\bx|\theta_0)p(\bx|\theta_0+\Delta)+
4qwp(\bx|\theta_0-\Delta)p(\bx|\theta_0+\Delta)]\Delta^2}{qp(\bx|\theta_0-\Delta)+
rp(\bx|\theta_0)+wp(\bx|\theta_0+\Delta)}.\nonumber
\end{eqnarray}
The lower bound on the MSE is then the sum of three integrals
\begin{eqnarray}
I_1&=&qr\Delta^2\cdot\int_{\reals^n}\frac{p(\bx|\theta_0-\Delta)p(\bx|\theta_0)\mbox{d}\bx}
{qp(\bx|\theta_0-\Delta)+rp(\bx|\theta_0)+wp(\bx|\theta_0+\Delta)}\\
I_2&=&rw\Delta^2\cdot\int_{\reals^n}\frac{p(\bx|\theta_0)p(\bx|\theta_0+\Delta)\mbox{d}\bx}
{qp(\bx|\theta_0-\Delta)+rp(\bx|\theta_0)+wp(\bx|\theta_0+\Delta)}\\
I_3&=&4qw\Delta^2\cdot\int_{\reals^n}\frac{p(\bx|\theta_0-\Delta)p(\bx|\theta_0+\Delta)\mbox{d}\bx}
{qp(\bx|\theta_0-\Delta)+rp(\bx|\theta_0)+wp(\bx|\theta_0+\Delta)}.
\end{eqnarray}
\\

\noindent
{\bf Example 2''.}
Revisiting again Example 2, for a given $t>0$, let us denote $\calC(t)=[0,t]^n$, and then
$p(\bx|\theta_0+i\Delta)=[\theta_0+i\Delta]^{-n}\cdot\calI\{\bx\in\calC(\theta_0+i\Delta)$,
$i=-1,0,1$.
Then,
\begin{eqnarray}
I_1&=&qr\Delta^2\cdot\int_{\calC(\theta_0-\Delta)}
\frac{(\theta_0-\Delta)^{-n}\theta_0^{-n}\mbox{d}\bx}
{q(\theta_0-\Delta)^{-n}+r\theta_0^{-n}+w(\theta_0+\Delta)^{-n}}\nonumber\\
&=&qr\Delta^2\cdot\frac{(\theta_0-\Delta)^n(\theta_0-\Delta)^{-n}\theta_0^{-n}}
{q(\theta_0-\Delta)^{-n}+r\theta_0^{-n}+w(\theta_0+\Delta)^{-n}}\nonumber\\
&=&\frac{qr\Delta^2\theta_0^{-n}}
{q(\theta_0-\Delta)^{-n}+r\theta_0^{-n}+w(\theta_0+\Delta)^{-n}}\nonumber\\
&=&\frac{qr\Delta^2}
{q[\theta_0/(\theta_0-\Delta)]^{n}+r+w[\theta_0/(\theta_0+\Delta)]^{n}}.
\end{eqnarray}
For $\Delta=\theta_0s/n$, we then have, for large $n$:
\begin{equation}
I_1\sim\frac{\theta_0^2}{n^2}\cdot\frac{qrs^2}{qe^s+r+we^{-s}}=
\frac{\theta_0^2}{n^2}\cdot\frac{qrs^2e^s}{qe^{2s}+re^s+w}.
\end{equation}
\begin{eqnarray}
I_2&=&rw\Delta^2\int_{\calC(\theta_0)}\frac{\theta_0^{-n}(\theta_0+\Delta)^{-n}\mbox{d}\bx}
{q(\theta_0-\Delta)^{-n}\calI\{\bx\in\calC(\theta_0-\Delta)\}+r\theta_0^{-n}+w(\theta_0+\Delta)^{-n}}\nonumber\\
&=&rw\Delta^2\cdot\frac{(\theta_0-\Delta)^n\theta_0^{-n}(\theta_0+\Delta)^{-n}}
{q(\theta_0-\Delta)^{-n}+r\theta_0^{-n}+w(\theta_0+\Delta)^{-n}}+\nonumber\\
& &rw\Delta^2\cdot\frac{[\theta_0^n-(\theta_0-\Delta)^n]\theta_0^{-n}(\theta_0+\Delta)^{-n}}
{r\theta_0^{-n}+w(\theta_0+\Delta)^{-n}}\nonumber\\
&=&rw\Delta^2\cdot\frac{(1-\Delta/\theta_0)^n}
{q[(1+\Delta/\theta_0)/(1-\Delta/\theta_0)]^{-n}+r(1+\Delta/\theta_0)^{n}+w}+\nonumber\\
& &rw\Delta^2\cdot\frac{[1-(1-\Delta/\theta_0)^n]}
{r(1+\Delta/\theta_0)^{n}+w}\nonumber\\
&\sim&\frac{\theta_0^2}{n^2}\cdot
rws^2\left(\frac{e^{-s}}{qe^{2s}+re^s+w}+\frac{1-e^{-s}}{re^s+w}\right).
\end{eqnarray}
Similarly,
\begin{equation}
I_3\sim \frac{\theta_0^2}{n^2}\cdot \frac{4qws^2}{qe^{2s}+re^s+w}.
\end{equation}
Thus,
\begin{eqnarray}
r(\theta_0)&\ge&\theta_0^2\cdot\sup_{s\ge 0}
\sup_{\{(q,r,w)\in\reals_+^3:~q+r+w=1\}} s^2\cdot\left[\frac{qre^s+4qw+rwe^{-s}}{qe^{2s}+re^s+w}+
\frac{rw(1-e^{-s})}{re^s+w}\right]\nonumber\\
&=&0.4624\theta_0^2,
\end{eqnarray}
w.r.t.\ $\zeta_n=n^2$,
which is a significant improvement of nearly 50\% over the previous bound,
$0.3102\theta_0^2$, that was obtained on the basis of two test points, let
alone the bound of Theorem 1, which was $0.2414\theta_0^2$.
This concludes Example 2''.\\

In general, the integrals $I_1$, $I_2$ and $I_3$ are not easy to calculate due
to the summations at the denominators of the integrands. One way to proceed is to
apply Lemma 1 to the sum of $k=3$ terms,
but this would introduce two additional parameters to be optimized.
Returning to the earlier shorthand notation of $a$, $b$, and $c$, a 
different approach to get rid of summations at the denominators, at the expense of some loss
of tightness, is the following: 
\begin{eqnarray}
\label{69}
\frac{ab+bc+4ac}{a+b+c}
&=&\frac{a(b+2c)}{a+b+c}+\frac{c(b+2a)}{a+b+c}\nonumber\\
&\ge&\frac{a(b+2c)}{a+b+2c}+\frac{c(b+2a)}{2a+b+c}\nonumber\\
&\ge&\frac{a(b+2c)}{2\max\{a,b+2c\}}+\frac{c(b+2a)}{2\max\{2a+b,c\}}\nonumber\\
&=&\frac{1}{2}\cdot\left(\min\{a,b+2c\}+\min\{2a+b,c\}\right)\nonumber\\
&\ge&\frac{1}{2}\cdot\left(\min\{a,b\}+\min\{b,c\}\right)\nonumber\\
&=&\frac{1}{2}\cdot\left[\min\{qp(\bx|\theta_0-\Delta),rp(\bx|\theta_0)\}+
\min\{rp(\bx|\theta_0),wp(\bx|\theta_0+\Delta)\}\right],
\end{eqnarray}
and so,
\begin{eqnarray}
R_n(\Theta)
&\ge&\frac{\Delta^2}{2}\bigg[\int_{\reals^n}
\min\{qp(\bx|\theta_0-\Delta),rp(\bx|\theta_0)\}\mbox{d}\bx+\nonumber\\
& &\int_{\reals^n}
\min\{rp(\bx|\theta_0),wp(\bx|\theta_0+\Delta)\}\mbox{d}\bx\bigg]\nonumber\\
&=&\frac{\Delta^2}{2}\bigg[(q+r)P_{\mbox{\tiny
e}}\left(\frac{q}{q+r},\theta_0-\Delta,\theta_0\right)+\nonumber\\
& &(r+w)P_{\mbox{\tiny
e}}\left(\frac{r}{r+w},\theta_0,\theta_0+\Delta\right)\bigg].
\end{eqnarray}
Note that by setting $w=0$, we recover the bound obtained by integrating
(\ref{37}), and therefore, by optimizing $w$, the resulting bound cannot be
worse. A slightly different (better, but more complicated) route in (\ref{69}) is 
to apply Lemma 1 with $k=2$ in the following manner:
\begin{eqnarray}
\frac{ab+bc+4ac}{a+b+c}
&=&\frac{a(b+2c)}{a+b+c}+\frac{c(b+2a)}{a+b+c}\nonumber\\
&\ge&\frac{a(b+2c)}{a+b+2c}+\frac{c(b+2a)}{2a+b+c}\nonumber\\
&=&\frac{a(b+2c)}{\min_u\max\{a/u,(b+2c)/(1-u)\}}+\frac{c(b+2a)}{\min_v\max\{(2a+b)/(1-v),c/v\}}\nonumber\\
&=&\max_u\min\{(1-u)a,u(b+2c)\}+\max_v\min\{v(2a+b),(1-v)c\}\nonumber\\
&\ge&\max_u\min\{(1-u)a,ub\}+\max_v\min\{vb,(1-v)c\}\nonumber\\
&=&\max_u\min\{(1-u)qp(\bx|\theta_0-\Delta),urp(\bx|\theta_0)\}+\nonumber\\
& &\max_v\min\{vrp(\bx|\theta_0),(1-v)wp(\bx|\theta_0+\Delta)\},
\end{eqnarray}
and so,
\begin{eqnarray}
R_n(\Theta)&\ge&\Delta^2\cdot\bigg[\max_u\int_{\reals^n}
\min\{(1-u)qp(\bx|\theta_0-\Delta),urp(\bx|\theta_0)\}\mbox{d}\bx+\nonumber\\
& &\max_v\int_{\reals^n}
\min\{vrp(\bx|\theta_0),(1-v)wp(\bx|\theta_0+\Delta)\}\mbox{d}\bx\bigg]\nonumber\\
&=&\Delta^2\cdot\bigg\{\max_u[(1-u)q+ur]\cdot P_{\mbox{\tiny
e}}\left(\frac{(1-u)q}{(1-u)q+ur},\theta_0-\Delta,\theta_0\right)+\nonumber\\
& &\max_v[vr+(1-v)w]\cdot P_{\mbox{\tiny
e}}\left(\frac{vr}{vr+(1-v)w},\theta_0,\theta_0+\Delta\right)\bigg\}.
\end{eqnarray}
We have just proved the following result:\\

\noindent
{\bf Theorem 4.} For the MSE case,
\begin{eqnarray}
R_n(\Theta)&\ge&\sup_{\theta_0,\Delta,q,r,w,(u,v)\in[0,1]^2}
\Delta^2\cdot\bigg\{\max_u[(1-u)q+ur]\cdot P_{\mbox{\tiny
e}}\left(\frac{(1-u)q}{(1-u)q+ur},\theta_0-\Delta,\theta_0\right)+\nonumber\\
& &\max_v[vr+(1-v)w]\cdot P_{\mbox{\tiny
e}}\left(\frac{vr}{vr+(1-v)w},\theta_0,\theta_0+\Delta\right)\bigg\}.
\end{eqnarray}
where the maximum over 
$q$, $r$, and $w$ is under the constraints that they are all non-negative and
sum to unity.\\

Revisiting Example 4, it is natural that both error probabilities would
correspond to a prior of $1/2$. This dictates the relations,
\begin{eqnarray}
u=\frac{q}{q+r}\\
v=\frac{w}{w+r},
\end{eqnarray}
and so, the bound becomes
\begin{eqnarray}
R_n(\Theta)&\ge&\max_{\Delta,q,w,r}
\Delta^2\left(\frac{2qr}{q+r}+\frac{2wr}{w+r}\right)\cdot
Q\left(\frac{\sqrt{n}\Delta}{2\sigma}\right)\nonumber\\
&=&\left[\max_{s\ge
0}\left(\frac{2s\sigma}{\sqrt{n}}\right)^2Q(s)\right]\cdot\max_{\{(q,r):~q\ge
0,~r\ge 0,~q+r\le
1\}}\left\{\frac{2qr}{q+r}+\frac{2r(1-q-r)}{1-q}\right\}\nonumber\\
&=&\frac{\sigma^2}{n}\cdot\max_{s\ge 0}\{4s^2Q(s)\}\cdot 0.6862\nonumber\\
&=&\frac{\sigma^2}{n}\cdot 0.6629\cdot 0.6862\nonumber\\
&=&\frac{0.4549\sigma^2}{n},
\end{eqnarray}
which is an improvement on the bound of Example 4, which was
$0.3314\sigma^2/n$. Similarly, in Example 5, for smooth signals, the
multiplicative constant improves from 0.1657 to 0.2274 and for the rectangular
pulse - from 0.1886 to 0.2588 (about 37\% improvement in all of them).

Revisiting Example 2,
we now have
\begin{eqnarray}
R_n(\Theta)&\ge&
\Delta^2\cdot[\max_u\min\{(1-u)q,ur\alpha\}+\max_v\min\{vr,(1-v)w\alpha\}]\nonumber\\
&=&\Delta^2\left(\frac{qr\alpha}{q+r\alpha}+\frac{wr\alpha}{r+w\alpha}\right)\nonumber\\
&=&\Delta^2\left(\frac{qr\alpha}{q+r\alpha}+\frac{r(1-q-r)\alpha}{r+(1-q-r)\alpha}\right)\nonumber\\
&=&\frac{s^2}{n^2}e^{-s}r\left(\frac{q}{q+re^{-s}}+\frac{1-q-r}{r+(1-q-r)e^{-s}}\right)\nonumber\\
&=&\frac{s^2}{n^2}r\left(\frac{q}{qe^s+r}+\frac{1-q-r}{re^s+(1-q-r)}\right),
\end{eqnarray}
whose maximum is $0.3909/n^2$, an improvement relative to the earlier bound of
$0.3102/n^2$.\\

\section*{Appendix A - Proof of eq.\ (\ref{approxrho})}
\renewcommand{\theequation}{A.\arabic{equation}}
    \setcounter{equation}{0}

Consider the Taylor series expansion,
\begin{equation}
s(t,\theta+\delta)=s(t,\theta)+\delta\cdot\dot{s}(t,\theta)+\frac{\delta^2}{2}\cdot
\ddot{s}(t,\theta)+o(\delta^2),
\end{equation}
where $\dot{s}(t,\theta)$ and $\ddot{s}(t,\theta)$ are the first two partial
derivatives of $s(t,\theta)$ w.r.t.\ $\theta$.
Correlating both sides with $s(t,\theta)$ yields
\begin{equation}
\int_0^Ts(t,\theta)s(t,\theta+\delta)\mbox{d}t=E+\delta\cdot\int_0^Ts(t,\theta)\dot{s}(t,\theta)\mbox{d}t
+\frac{\delta^2}{2}\cdot\int_0^Ts(t,\theta)\ddot{s}(t,\theta)\mbox{d}t+o(\delta^2).
\end{equation}
Now,
\begin{equation}
\int_0^Ts(t,\theta)\dot{s}(t,\theta)\mbox{d}t=\frac{1}{2}\cdot\frac{\partial}{\partial\theta}\left\{
\int_0^Ts^2(t,\theta)\mbox{d}t\right\}=\frac{1}{2}\cdot\frac{\partial
E}{\partial\theta}=0,
\end{equation}
since the energy $E$ is assumed independent of $\theta$. Also, since
$\frac{\partial^2E}{\partial\theta^2}=0$, we have
\begin{eqnarray}
0&=&\frac{\partial^2}{\partial\theta^2}\left\{\int_0^Ts^2(t,\theta)\mbox{d}t\right\}\nonumber\\
&=&\frac{\partial}{\partial\theta}\left\{2\cdot\int_0^Ts(t,\theta)\dot{s}(t,\theta)\mbox{d}t\right\}\nonumber\\
&=&2\cdot\int_0^T\dot{s}^2(t,\theta)\mbox{d}t+2\cdot\int_0^Ts(t,\theta)\ddot{s}(t,\theta)\mbox{d}t,
\end{eqnarray}
which yields
\begin{equation}
\int_0^Ts(t,\theta)\ddot{s}(t,\theta)\mbox{d}t=-\int_0^T\dot{s}^2(t,\theta)\mbox{d}t,
\end{equation}
and so,
\begin{equation}
\int_0^Ts(t,\theta)s(t,\theta+\delta)\mbox{d}t=E-\frac{\delta^2}{2}\int_0^TT\dot{s}^2(t,\theta)\mbox{d}t+o(\delta^2),
\end{equation}
and so,
\begin{equation}
\varrho(\theta,\theta+\Delta)=1-\frac{\delta^2}{2E}\int_0^TT\dot{s}^2(t,\theta)\mbox{d}t+o(\delta^2).
\end{equation}

\section*{Appendix B -- Proof of Lemma 1}
\renewcommand{\theequation}{B.\arabic{equation}}
    \setcounter{equation}{0}

First, observe that
\begin{equation}
\sum_{i=1}^ka_i=\sum_{i=1}r_i\cdot\frac{a_i}{r_i}\le\max\left\{\frac{a_1}{r_1},\ldots,\frac{a_k}{r_k}\right\},
\end{equation}
and since the inequality $\sum_{i=1}^ka_i\le \max\{a_1/r_1,\ldots,a_k/r_k\}$
applies to all $\br\in\calS$, it also applies to the infimum 
of $\max\{a_1/r_1,\ldots,a_k/r_k\}$ over $\calS$, thus establishing 
the inequality ``$\le$'' between the two sides. To establish the
``$\ge$'' inequality, define $\br^*\in\calS$ to be the vector whose components
are given by $r_i^*=a_i/\sum_{j=1}^ka_j$. Then,
\begin{equation}
\inf_{(r_1,\ldots,r_k)\in\calS}\max\left\{\frac{a_1}{r_1},\ldots,\frac{a_k}{r_k}\right\}\le
\max\left\{\frac{a_1}{r_1^*},\ldots,\frac{a_k}{r_k^*}\right\}=\sum_{i=1}^ka_i.
\end{equation}
This completes the proof of Lemma 1.

\end{document}